\documentclass[a4paper,12pt]{article}
\usepackage{hyperref}
\usepackage[english]{babel}
\usepackage[T1]{fontenc}
\usepackage[utf8]{inputenc}
\usepackage{amsfonts}
\usepackage{amssymb}
\usepackage{amsmath}
\usepackage{graphicx}
\usepackage{cite}
\usepackage{epsfig}
\usepackage{psfrag}
\usepackage{tikz}
\usetikzlibrary{matrix,arrows,decorations.pathmorphing,positioning} 
\usepackage{url}
\usepackage{amsthm}
\usepackage{mathpartir}
\usepackage{enumerate}
\usepackage{tikz-cd}

\newtheorem{proposition}{Proposition}[section]
\newtheorem{definition}{Definition}[section] 
  
\newtheorem{corollary}{Corollary}[section] 
\newtheorem{remark}{Remark}[section]  

\date{}

\title{\bf Pregeometric Spaces from Wolfram Model Rewriting Systems as   Homotopy Types}

\author{Xerxes D. Arsiwalla$^{1, 3, }$\footnote{Corresponding Author: \url{x.d.arsiwalla@gmail.com}}  {}  and  Jonathan Gorard$^{2, 3, }$\footnote{\url{jg865@cam.ac.uk}}\\  
{}  \\
{\it \small $^{1}$Pompeu Fabra University, Barcelona, Spain}\\ 
{\it \small $^{2}$University of Cambridge, Cambridge, United Kingdom}\\
{\it \small $^{3}$Wolfram Research, USA}    
}

\begin{document}
\maketitle

\begin{abstract}
How do spaces in physics emerge from pregeometric discrete building blocks governed by computational rules?  To address this question we investigate  non-deterministic rewriting systems (so-called multiway systems) of the Wolfram model. We express these rewriting systems as homotopy types. Using this new formulation of the model, we motivate  how spatial structures are  functorially inherited from pregeometric type-theoretic constructions. We show how higher homotopies are constructed from rewriting rules. These correspond to morphisms of an $n$-fold category. Subsequently, the $n \to \infty$ limit of the Wolfram model rulial multiway system is identified as an $\infty$-groupoid, with the latter being relevant  in light of Grothendieck's homotopy hypothesis.  We then go on to show how this construction extends to the classifying space of rulial multiway systems, which forms a multiverse of multiway systems and carries the formal structure of an ${\left(\infty, 1\right)}$-topos. This correspondence to higher categorical structures potentially offers a new way to understand how the kinds of spaces relevant to physics result from pregeometric  combinatorial models.  The key issue we have addressed in this work is to relate abstract non-deterministic rewriting systems to higher homotopy spaces.  A  consequence of constructing spaces and geometry synthetically is that it removes the need to make ad hoc assumptions about  geometric attributes of a model such as an a priori background  or any pre-assigned geometric data. Instead, geometry is inherited functorially from higher  structures.  This can be particularly useful for formally justifying different choices of underlying spacetime discretization schemes adopted by various models of quantum gravity.   We conclude with comments on how our   framework of higher category-theoretic combinatorial constructions closely corroborates  with  other  approaches investigating higher categorical structures relevant to the foundations of physics.  

\end{abstract}

\clearpage

\tableofcontents

\clearpage

\section{Introduction}


One of the major challenges facing the foundations of physics is the problem of reconciling quantum theory with general relativity, arguably, the two  hallmarks of contemporary fundamental physics. There are reasons to believe that such a reconciliation will markedly shape our understanding of several open  problems in theoretical physics, including, the origin of our universe, the origin of matter  and  the fundamental forces, quantum properties of black holes, the nature of space and time; among others. In a sense, a key question underlying many of these issues is the following: What are the fundamental building blocks of the universe? Any theory attempting to reconcile quantum mechanics with general relativity will, at the very least, have to take a stance on the nature of these building blocks (those referring to spacetime and/or matter) and then proceed thereon. While it certainly seems there is a consensus to the view that space, time and matter are fundamentally discrete \cite{isham1995structural}, \cite{gibbs1995small};  specific proposals concerning the nature of this discretization (sometimes traded for quantization) and consequently the underlying mathematical structure one needs to start with, wildly differ. Notable examples of such efforts include (i) theories of quantum gravity such as loop quantum gravity \cite{rovelli2008loop},  manifestations of string theory and M-theory as quantized fluctuations about a  background    \cite{becker2006string}, causal dynamical triangulations \cite{loll2019quantum}, causal set theory  \cite{dowker2006causal},  \cite{rideout2009emergence};  (ii) recent  approaches to quantization invoking grand unification within supersymmetric  field  theories  \cite{altarelli20015},  F-theory  \cite{beasley2009gutsI}, \cite{beasley2009gutsII}; and extended spin foam models \cite{lisi2010unification};  (iii) models of emergent spacetime such as spacetime from entanglement  \cite{swingle2018spacetime},  AdS/CFT, gauge-gravity duality and related holographic models \cite{aharony2000large}, \cite{arsiwallasupersymmetric},  \cite{arsiwalla2006phase},     \cite{arsiwalla2011degenerate},  emergent gravity \cite{verlinde2017emergent},  energetic causal sets \cite{cortes2014universe}; (iv) proposals of pregeometric physics, most notably those championed by Wheeler \cite{misner1973gravitation}, \cite{wheeler1980pregeometry}

On one hand, it is interesting to note that many of the theories cited above (if  not all of them) come equipped with a fair share of a priori geometric structure.   On the other hand, a true pregeometric  description ought to be one from which all geometric features of the physical universe ought to be derived (as extensively argued for in \cite{meschini2005geometry}).  In this view, it is the precursors of geometry (and to a certain extent, even topology) that  make up the abstract building blocks of the universe. The term `pregeometry' was first coined by John Wheeler as a non-geometric approach that ought to   encompass  any underlying explanation of spacetime or quantum gravity  \cite{misner1973gravitation}, \cite{wheeler1980pregeometry}.  The argument  can be stated as follows:  given that quantum mechanics permits metrics to fluctuate,   merging gravity with quantum mechanics, at the very least, requires a set of more fundamental rules regarding connectivity of spacetime that are independent of topology and dimensionality. While contemporary formulations of physics often come pre-defined with a priori geometry, formulations based on  pregeometric structures may  allow one to work with deeper underlying rules of physics that are not dependent on the usual structural assumptions about the properties of space and time.  In this work, we undertake the problem of formalizing pregeometric spaces using the framework of homotopy type theory and higher categories. This will necessitate building conceptual  bridges between algebraic logic, proof theory, formal models of computation and physics (or rather pregeometric physics). These connections may potentially offer us new insights underlying the formal foundations of quantum process theories such as categorical quantum mechanics as well as discrete models of  spacetime.

This article is an attempt at addressing two key objectives:   Firstly, we seek to provide a mathematical framework for generating formal  constructions of  pregeometric structures, which, in certain to-be-specified limits, would yield candidate spaces within which topological or geometric structures relevant to physics may be modeled.   
The other objective of this work is to provide a new formal foundation for the Wolfram model  \cite{Wolfram2020},  \cite{Wolfram2002a} (see also \cite{Gorard2020}, \cite{Gorard2020a}, \cite{Gorard2020c}, \cite{arsiwalla2021homotopy}) using homotopy types.  Both these issues are, of course, related.  The Wolfram model originates from the idea that  the building blocks of the universe are fundamentally discrete entities (and their relations), governed by computational rules. Structures relevant to physics are subsequently constructed in this model from the combinatorics of non-deterministic rewriting systems. Using our new formulation of this model using  homotopy types, we show how spatial structures can be realized  from suitable limits of pregeometric homotopy types in an $( \infty, 1 )$ category.  The Wolfram model formalized as a homotopical type theory thus provides a new perspective on Wheeler's original intuition of physics from pregeometry.  We reckon this framework of higher categorical pregeometric structures may also be relevant to several other approaches addressing questions at the foundations of physics; for instance, the higher pre-quantum geometry program  \cite{schreiber2016higher}, higher gauge theories on cohesive toposes  \cite{Schreiber2012},  \cite{Schreiber2013a}, extended topological and axiomatic quantum field theories \cite{lurie2009classification}, \cite{schreiber2009aqft},  categorical quantum mechanics \cite{abramsky2009categorical},  \cite{coecke2018picturing},  among others.



A brief word about the Wolfram model: The Wolfram model  \cite{Wolfram2020},  \cite{Wolfram2002a}  is an explicitly computational framework based on non-deterministic rewriting systems. This model attempts to provide a description of fundamental physical phenomenon starting from rewriting systems.  It seeks to capture ways in which simple rewriting rules may be composed, so as to yield more complex structures (suggestively titled ``universes")  that admit certain emergent ``laws of physics''. This general approach of seeking a constructivist formalization of  physical  structures has strong parallels in recent work in the foundations of both physics and mathematics, for instance in relation to synthetic geometry and cohesive homotopy type theory \cite{Schreiber2012}, \cite{Schreiber2013a},  \cite{Shulman2016},  \cite{Ahrens2021}.   The general framework of the Wolfram model eschews the rigid continuum description of spacetime in terms of Lorentzian manifolds in favor of a more rudimentary description of the combinatorics and relational properties of intrinsic causality between events,  which, in turn, makes manifest the computational architecture that underlies contemporary physics. The archetypical structures that appear within this framework are so called `multiway systems' - non-deterministic abstract rewriting systems equipped with a notion of causal structure. The `multiway' moniker  designates the fact that all permissible applications of  rewriting rules are instantiated in all possible orderings, leading to multiple chains of rewriting terms that are partially-ordered by causality. Depending upon the precise interpretation of these terms or multiway states, Wolfram model multiway systems may be realized as rewriting systems over graphs, hypergraphs,    character strings, ZX diagrams \cite{Gorard2020c}, \cite{Gorard2021a} (as formalized by Coecke and Duncan \cite{Coecke2008}, \cite{Coecke2009a}), string diagrams  \cite{Gorard2021b}  (as formalized by Joyal and Street \cite{Joyal1991}), etc.

In this work, we further investigate non-deterministic rewriting systems. For our purposes, the states of these systems can be purely abstract entities. That will allow for studying algebraic and compositional properties of these systems at the highest level of generality. Consequently,  results thus obtained, may be easily specialized to relevant classes of multiway rewriting systems. Specifically, we show how Wolfram model multiway systems can be enhanced with higher homotopies. These homotopies are induced by  inclusion of   `higher-order' rules taken from a so-called `rulial space', which designates the space of all possible rewriting rules of a given signature (for the case of graphs or hypergraphs, this space would comprise a monoidal category of cospans). Although a generic  rewriting system may be thought of as simply being an $F$-coalgebra for the power set functor  ${\mathcal{P}}$  (which, in the case of hypergraph rewriting, is further equipped with a natural symmetric monoidal structure), in earlier work  \cite{arsiwalla2021homotopy},  we showed that a rulial multiway rewriting system equipped with homotopies up to order $n$ can elegantly be formalized as an $n$-fold category. Furthermore, upon including inverse morphisms (via the inclusion of invertible rewriting relations), the $n \to \infty$ limit of the rulial multiway system yields an ${\infty}$-groupoid, which inherits the structure of a formal homotopy space via Grothendieck's homotopy hypothesis. Under certain conditions, these structures inherit non-trivial topologies, and thus, establish formal connections between the combinatorial framework of Wolfram model rewriting systems and topological spaces (with additional of a cohesive structure, one also obtains  geometric spaces). Hence, starting from purely pregeometric rewriting constructs, we attempt to chart out a formal route to generating spaces, upon which the laws of physics may eventually be realized.

Moreover, the aforementioned connections of non-deterministic rewriting systems to higher homotopy theory may be incorporated within the broader context of Homotopy Type Theory (HoTT). We go on to establish the correspondence between Wolfram model multiway systems and higher category theory. We outline how combinatorial structures in the Wolfram model, such as  multiway rewriting systems, rulial multiway systems, branchial graphs, causal networks, etc., can be expressed as homotopy types.  In particular, we explicitly represent multiway rewriting systems using type constructors and discuss how Wolfram model constructions as homotopy types can be internalized within a suitable $\infty$-topos, within which, spaces and constructions relevant to physics can potentially be realized.  Furthermore, we show how this construction extends to the classifying space of rulial multiway systems, which forms a multiverse of multiway systems and carries the formal structure of an ${\left(\infty, 1\right)}$-topos.  This correspondence to spaces and higher structures  offers a formal understanding of how  spatial structures relevant to theories  of physics can be shown to originate  from abstract combinatorial constructions   \cite{Arsiwalla2020}, \cite{arsiwalla2021homotopy},  \cite{Shulman2017}, \cite{Baez2006}.  A related program seeking to formalize quantum field theories from cohesive $\infty$-toposes also borrows heavily from such a synthetic approach to geometry based on higher categories   \cite{Schreiber2012},  \cite{Schreiber2013a}.

The outline of this paper is as follows: In section 2, we introduce the preliminaries of the Wolfram model as an abstract rewriting system, its 1-categorical formulation and an overview of other related works. In section 3, we provide the detailed type-theoretic construction of multiway systems, and also show that mathematical objects as groups, categories and groupoids can equivalently be expressed as type constructions or Wolfram model constructions. In section 4, we present an algorithmic framework for constructing higher homotopical structures in multiway rewriting systems. In earlier work,   we showed that these realize $n$-fold categories, and ${\infty}$-groupoids in limiting cases. Here, we extend that to a multiverse of multiway systems, which has the formal structure of an ${\left(\infty, 1\right)}$-topos. In section 5, we outline how the above type-theoretic structures can be internalized within an $\infty$-topos, and discuss connections to the synthetic geometry program. In section 6, we discuss potential applications of our approach to physics; namely, how spatial structures relevant to physics may arise from rewriting systems, connections to categorical quantum mechanics and topological field theories, the essential role of the observer in a constructivist paradigm, and a homotopical interpretation of graph and hypergraph limits for discrete spacetime models. Finally, in section 7, we conclude with final remarks and future directions.


%
%
%

\clearpage

\section{Preliminaries of the Wolfram Model }

\subsection{The Wolfram Model as an Abstract Rewriting System}

\begin{figure}[ht]
\centering
\includegraphics[width=0.2\textwidth]{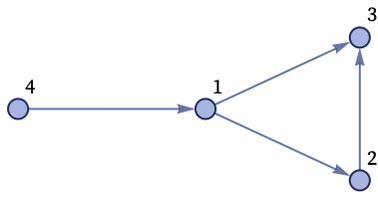}\hspace{0.25\textwidth}
\includegraphics[width=0.2\textwidth]{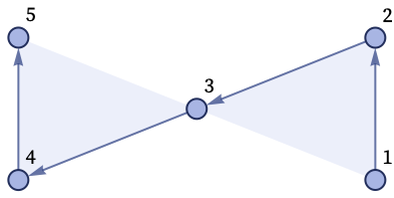}
\caption{Spatial hypergraphs corresponding to finite collections of ordered relations between elements, namely ${\left\lbrace \left\lbrace 1, 2 \right\rbrace, \left\lbrace 1, 3 \right\rbrace, \left\lbrace 2, 3 \right\rbrace, \left\lbrace 4, 1 \right\rbrace \right\rbrace}$ and ${\left\lbrace \left\lbrace 1, 2, 3 \right\rbrace, \left\lbrace 3, 4, 5 \right\rbrace \right\rbrace}$, respectively.}
\label{fig:Figure1}
\end{figure}

We begin with a preliminary description of the Wolfram model in terms of diagrammatic rewriting rules acting on  hypergraphs. This model is a discrete spacetime formalism. It posits that structures such as continuous spacetime geometries and Hilbert spaces may potentially emerge from large-scale limits of the underlying discrete structures. The discrete structures of this model are hypergraphs, causal networks, branchial graphs and multiway graphs. Furthermore, the evolution of these structures is dictated by graph, hypergraph or string rewriting rules.   
A Wolfram model hypergraph can be represented abstractly as finite collections of ordered relations (i.e. hyperedges) between elements (i.e. hypernodes), as defined below and shown in Figure \ref{fig:Figure1}:

\begin{definition}
A  Wolfram model hypergraph, denoted ${H = \left( V, E \right)}$, is a finite collection of hyperedges (ordered):
\begin{equation}
E \subset \mathcal{P} \left( V \right) \setminus \left\lbrace \emptyset \right\rbrace,
\end{equation}
where ${\mathcal{P} \left( V \right)}$ denotes the power set of $V$ and $E$ is an ordered collection of nodes.
\end{definition}

\begin{figure}[ht]
\centering
\includegraphics[width=0.395\textwidth]{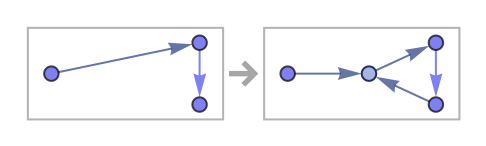}
\caption{A hypergraph transformation rule corresponding to the set substitution system ${\left\lbrace \left\lbrace x, y \right\rbrace, \left\lbrace y, z \right\rbrace \right\rbrace \to \left\lbrace \left\lbrace w, y \right\rbrace, \left\lbrace y, z \right\rbrace, \left\lbrace z, w \right\rbrace, \left\lbrace x, w \right\rbrace \right\rbrace}$.}
\label{fig:Figure2}
\end{figure}

One can then define the dynamics of a Wolfram model system in terms of hypergraph rewriting rules as follows:

\begin{definition}
An `update rule', denoted $R$, for a spatial hypergraph ${H = \left( V, E \right)}$ is an abstract rewriting rule of the form ${H_1 \to H_2}$, in which a subhypergraph matching pattern ${H_1}$ is replaced by a distinct subhypergraph matching pattern ${H_2}$.
\end{definition}
Each such rewriting rule is formally equivalent to a set substitution system (one in which a subset of ordered relations matching a particular pattern is replaced with a distinct subset of ordered relations matching a particular pattern), as shown in Figure \ref{fig:Figure2}.

Note that, in general, the order in which to apply the transformation rules is not well-defined; in the simplest case, we could simply apply the rule to every possible matching (and non-overlapping) subhypergraph, as illustrated in Figures \ref{fig:Figure3} and \ref{fig:Figure4}. However, even in this simplified case, the initial choice of the subhypergraph to which to apply the first transformation is still ambiguous, and different such choices will in general yield non-isomorphic sequences of hypergraphs in the evolution.

\begin{figure}[ht]
\centering
\includegraphics[width=0.595\textwidth]{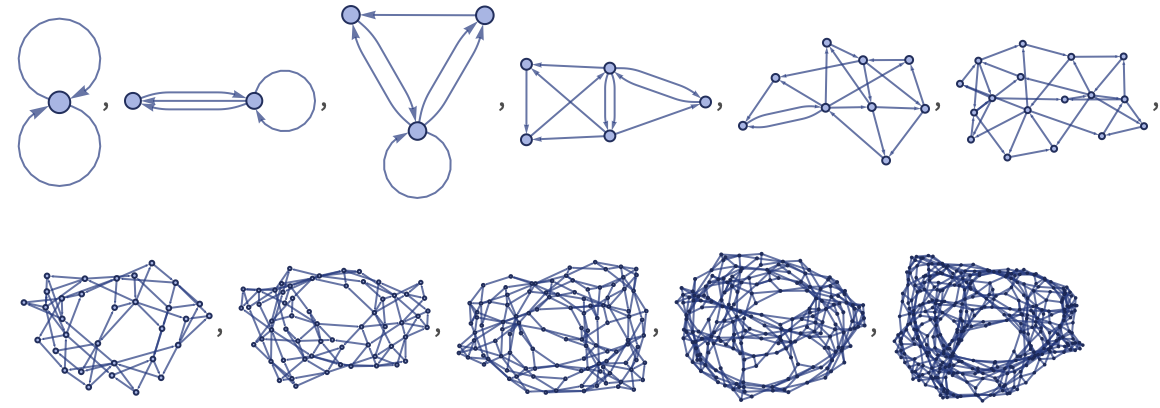}
\caption{The results of the first 10 steps in the evolution history of the set substitution system ${\left\lbrace \left\lbrace x, y \right\rbrace, \left\lbrace y, z \right\rbrace \right\rbrace \to \left\lbrace \left\lbrace w, y \right\rbrace, \left\lbrace y, z \right\rbrace, \left\lbrace z, w \right\rbrace, \left\lbrace x, w \right\rbrace \right\rbrace}$, starting from a double self-loop initial condition.}
\label{fig:Figure3}
\end{figure}

\begin{figure}[ht]
\centering
\includegraphics[width=0.3\textwidth]{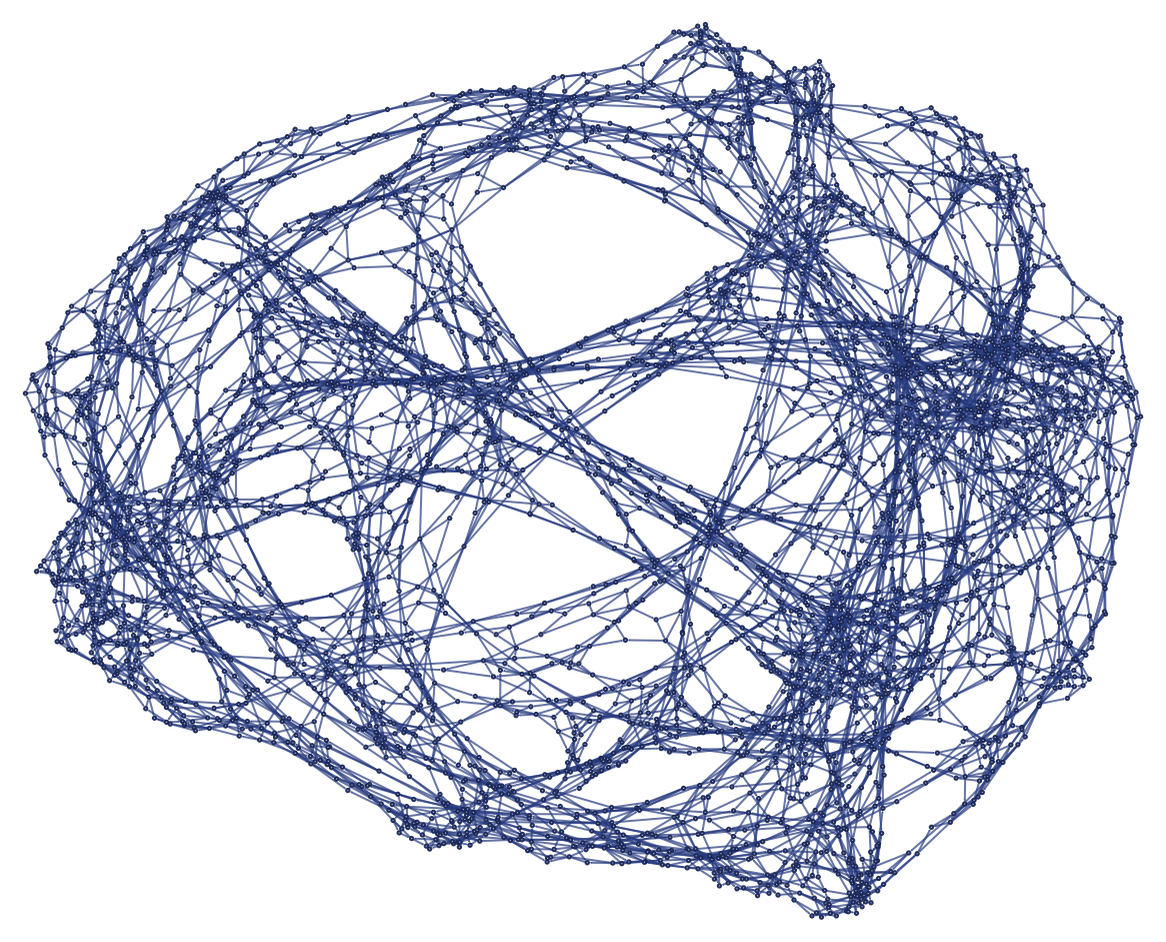}
\caption{The result after 14 steps of evolution of the set substitution system ${\left\lbrace \left\lbrace x, y \right\rbrace, \left\lbrace y, z \right\rbrace \right\rbrace \to \left\lbrace \left\lbrace w, y \right\rbrace, \left\lbrace y, z \right\rbrace, \left\lbrace z, w \right\rbrace, \left\lbrace x, w \right\rbrace \right\rbrace}$, starting from a double self-loop initial condition.}
\label{fig:Figure4}
\end{figure}

Therefore, the evolution of any given spatial hypergraph will, generically, be   non-deterministic, due to this lack of any canonical updating order; we can parametrize this non-determinism by treating the Wolfram model as an abstract rewriting system \cite{Baader1998}, \cite{Bezem2003}.

\begin{definition}
An `abstract rewriting system' (or `ARS') is a set, denoted $A$ (with each element known as an `object'), equipped with some binary relation, denoted ${\to}$, known as the `rewriting relation'.
\end{definition}

A concrete way of representing the abstract rewriting structure of a Wolfram model system is through the use of a general combinatorial structure known as a \textit{multiway rewriting system}. Combinatorially, a multiway rewriting system or a multiway system is simply a directed, acyclic graph of states, determined by abstract rewriting rules that inductively generate a (potentially infinite) \textit{multiway evolution graph}, together with a partial order on the rewriting rule applications, determined by their causal structure. This is straightforward to formalize in terms of \textit{abstract rewriting systems} \cite{Baader1998}, \cite{Bezem2003} in which the underlying rewriting relation is not (necessarily) confluent \cite{Dershowitz1990}, \cite{Huet1980}.

\begin{definition}
A `multiway evolution graph', denoted ${G_{multiway}}$, is a directed, acyclic graph corresponding to the evolution of a (generically non-confluent) abstract rewriting system ${\left( A, \to \right)}$, in which the set of vertices corresponds to the set of objects ${V \left( G_{multiway} \right)}$, and in which the directed edge ${a \to b}$ exists in ${E \left( G_{multiway} \right)}$ if and only if there exists an application of the rewriting relation that transforms object $a$ to object $b$.
\end{definition}
Hence, directed edges will connect vertices $a$ and $b$ in ${G_{multiway}}$ if and only if ${a \to b}$ in the underlying rewriting system, and a directed path will connect vertices $a$ and $b$ if and only if ${a \to^{*} b}$, where ${\to^{*}}$ denotes the reflexive transitive closure of ${\to}$, i.e if and only if there exists a finite rewriting sequence of the form:

\begin{equation}
a \to a^{\prime} \to a^{\prime \prime} \to \cdots \to b^{\prime} \to b
\end{equation}

Thus, the evolution of a generic Wolfram model system will correspond to a multiway evolution graph, within which the ``standard'' updating order shown above will correspond to a single path, as illustrated in Figure \ref{fig:Figure5}.  

\begin{figure}[ht]
\centering
\includegraphics[width=0.695\textwidth]{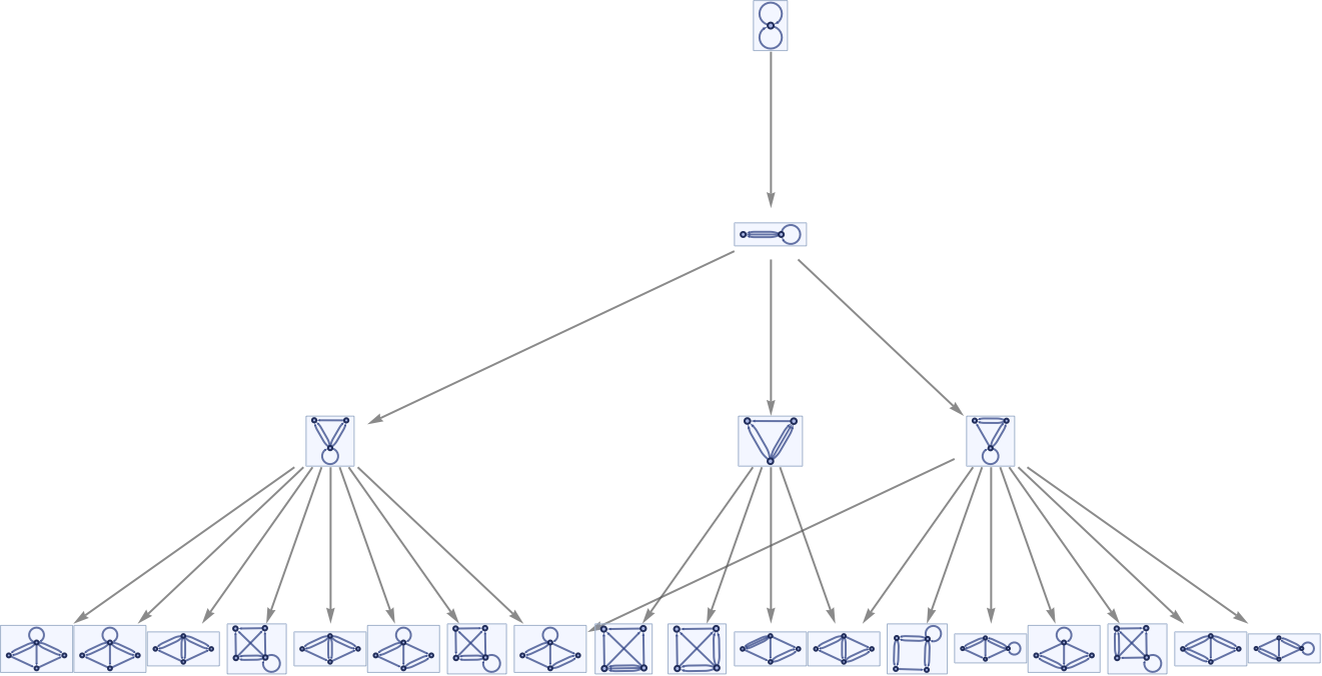}
\caption{The multiway evolution graph corresponding to the first 3 steps in the non-deterministic evolution history of the hypergraph substitution rule ${\left\lbrace \left\lbrace x, y \right\rbrace, \left\lbrace y, z \right\rbrace \right\rbrace \to \left\lbrace \left\lbrace w, y \right\rbrace, \left\lbrace y, z \right\rbrace, \left\lbrace z, w \right\rbrace, \left\lbrace x, w \right\rbrace \right\rbrace}$, starting from a simple double self-loop initial condition ${\left\lbrace \left\lbrace 0, 0 \right\rbrace, \left\lbrace 0, 0 \right\rbrace \right\rbrace}$.}
\label{fig:Figure5}
\end{figure}

Furthermore, in the case that future applications of transformation rules in the Wolfram model may have dependencies upon prior rule applications, such that updating event $b$ could only have been applied if event $a$ had previously been applied; such dependencies may be captured by means of a causal network:

\begin{definition}
A `causal network', denoted ${G_{causal}}$, is a directed, acyclic graph in which every vertex corresponds to an application of an update rule (i.e. an update `event'), and in which the directed edge ${a \to b}$ exists if and only if:

\begin{equation}
\mathrm{In} \left( b \right) \cap \mathrm{Out} \left( a \right) \neq \emptyset,
\end{equation}
i.e. the input for event $b$ makes use of hyperedges that were produced by the output of event $a$.
\end{definition}
In the context of the Wolfram model, the transitive reduction of a causal network is presumed to correspond to the Hasse diagram of the causal partial order for some discretized approximation to spacetime.

\begin{figure}[ht]
\centering
\includegraphics[width=0.395\textwidth]{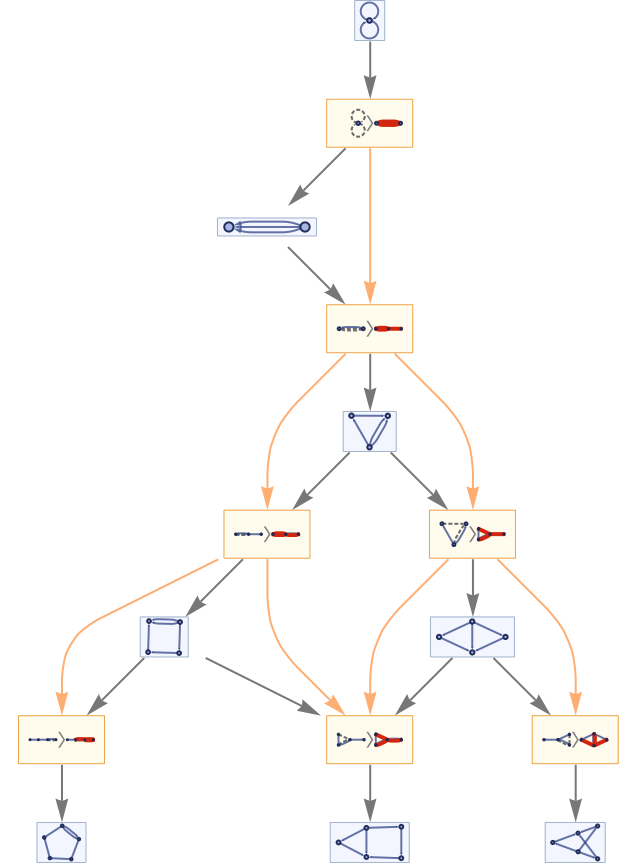}
\caption{The multiway evolution causal graph (with evolution edges shown in gray, and causal edges shown in orange) for the set substitution system ${\left\lbrace \left\lbrace x, y \right\rbrace, \left\lbrace z, y \right\rbrace \right\rbrace \to \left\lbrace \left\lbrace x, w \right\rbrace, \left\lbrace y, w \right\rbrace, \left\lbrace z, w \right\rbrace \right\rbrace}$, illustrating trivial causal invariance.}
\label{fig:Figure35}
\end{figure}

An example of a multiway evolution causal graph (in which updating events are shown in yellow, state vertices are shown in blue, evolution edges are shown in gray and causal edges are shown in orange) for a system exhibiting trivial causal invariance is featured in Figure \ref{fig:Figure35}.

The notion of confluence in abstract rewriting theory is deeply related to (and, indeed, is a necessary but not sufficient condition for) the criterion of \textit{causal invariance} in multiway evolution:

\begin{definition}
`Causal Invariance' is defined as a property of multiway systems whereby all possible evolution paths yield causal networks that are (eventually) isomorphic as directed acyclic graphs. 
\end{definition}

Since the notion of confluence in the theory of abstract rewriting systems is a necessary (though not sufficient) condition for causal invariance, it follows that whenever causal invariance exists, every branch in the multiway evolution graph must eventually merge. For the particular case of a terminating (strongly normalizing) rewriting system, causal invariance therefore implies that all evolution paths yield the same eventual state. The physical significance of causal invariance is believed to be closely tied to spacetime symmetries in the continuum limit  \cite{Wolfram2020}.

For later reference, the following definitions will also be useful:

\begin{definition}
 A  `\textit{Branchial Graph}'  is a  graph whose vertex set is the set of states in a particular layer (or \textit{slice}) of the multiway evolution graph, and in which states are connected by directed edges if and only if they share a common ancestor in the evolution graph. Otherwise known as a \textit{branchlike hypersurface}, by analogy to spacelike hypersurfaces in causal networks.  \end{definition}

In the Wolfram model, branchial graphs are used to indicate instantaneous superpositions between pure states.

\begin{definition}
`\textit{Branchial Space}'  encompasses the corresponding spatial  structure defined by a branchial graph, much like how physical space is the spatial structure defined by a hypergraph. In this way, branchial space has the same relationship to the multiway evolution graph as physical space has to an ordinary causal network. 
\end{definition}

\begin{definition}
A `\textit{Foliation}'  in the Wolfram model  is a method for defining a universal time function over the vertices of a directed acyclic graph (i.e., a function mapping vertices to integers), in such a way that the level sets of that function, known as \textit{slices}, cover the entire graph without intersecting. Foliations of a causal network yield successive configurations of hypergraphs, representing spacelike hypersurfaces. Foliations of a multiway evolution graph yield successive configurations of branchial graphs, representing instantaneous superpositions between pure states (\textit{branchlike hypersurfaces}).
\end{definition}

Software for reproducing all the above Wolfram model objects can be found here:  \url{https://www.wolframcloud.com/obj/wolframphysics/Tools/guide-page}

\subsection{Categorical Formulation of Wolfram Model (Hyper)graph Rewriting}


What we have presented above was a purely combinatorial description of the Wolfram model in terms of graph/hypergraph transformation rules as elementary operations on sets. We now provide a (1-)category description of the same class of hypergraph transformations. This formalization provides a compositional description in terms of $F$-coalgebras and a rewriting description  in terms of  double-pushout (DPO) rewriting over (selective) adhesive categories; the latter, follows from the construction in   \cite{Dixon2013}.  The contents of this subsection serve merely as a preliminary introduction to known constructions of Wolfram model multiway systems; the informed reader can safely proceed to the following sections. 


In the most general case, in which the rewriting relation ${\to}$ is treated as an indexed union of sub-relations, i.e. ${\to = \to_1 \cup \to_2 \cup \dots}$, with label set ${\Lambda}$ for the indices, and in which the set of objects $A$ is arbitrary (i.e. its elements could represent graphs, hypergraphs, string diagrams, terms, character strings, etc.), the resulting labeled abstract rewriting system ${\left( A, \Lambda, \to \right)}$ permits an elegant compositional description in terms of $F$-coalgebras. Specifically, note that the system ${\left( A, \Lambda, \to \right)}$ is now simply a bijective function from $A$ to a subset of the power set of $A$ indexed by ${\Lambda}$, i.e. ${\mathcal{P} \left( \Lambda \times A \right)}$:

\begin{equation}
p \mapsto \left\lbrace \left( \alpha, q \right) \in \Lambda \times A : p \to^{\alpha} q \right\rbrace.
\end{equation}
Recall now that an $F$-coalgebra for an endofunctor ${F : \mathbf{C} \to \mathbf{C}}$ consists of an object $A$ in ${\mathrm{ob} \left( \mathbf{C} \right)}$ equipped with a morphism ${\alpha : A \to F A}$ in ${\mathrm{hom} \left( \mathbf{C} \right)}$, hence denoted ${\left( A, \alpha \right)}$. Thus, since the power set construction on ${\mathbf{Set}}$ is a covariant endofunctor ${\mathcal{P} : \mathbf{Set} \to \mathbf{Set}}$, we see that the abstract rewriting system ${\left( A, \to \right)}$ consists of an object $A$ equipped with an additional morphism of ${\mathbf{Set}}$, namely the rewriting relation ${\to}$:

\begin{equation}
\to : A \to \mathcal{P} A.
\end{equation}


Furthermore, for a large class of cases, such as graph/hypergraph rewriting, open graph rewriting, string diagram rewriting (including the ZX-calculus), etc., the rewriting relation ${\to}$ itself may be further specified via  double-pushout (DPO) rewriting. In \cite{Dixon2013}  this was formalized  for the case of open graphs using a selective adhesive category of open spans, and the same formalism was subsequently applied to both ZX diagrams and arbitrary hypergraph rewriting systems as special cases  \cite{Gorard2020c}, \cite{Gorard2021a},  \cite{Gorard2021b}. And more recently, the mathematical foundations for both, double-pushout string diagram rewriting modulo Frobenius structures and double-pushout hypergraph rewriting have been discussed in  \cite{bonchi2020stringI},  \cite{bonchi2021stringII}  and   \cite{bonchi2021stringIII}. 

Within this framework, one can formally specify rewriting rules as spans of monomorphisms ${\rho}$ of the form:
\begin{equation}
\rho = \left( l : K \to L, r : K \to R \right),
\end{equation}
with the left- and right-hand-sides of the rule specified by objects $L$ and $R$, respectively, and with object $K$ designating the interface graph. A \textit{match} for the rule ${\rho}$ within an object $G$ is simply a morphism ${m : L \to G}$, with the rule being applicable if and only if there exists a pair of pushout diagrams of the form \cite{Ehrig1973}, \cite{Habel2001}:

\begin{equation}
\begin{tikzcd}
L \arrow[d, "m"] & K \arrow[l, "l"] \arrow[d, "k"] \arrow[r, "r"] & R \arrow[d, "n"]\\
G & D \arrow[l, "f"] \arrow[r, "g"] & H
\end{tikzcd}.
\end{equation}
These rules are defined and applied in the context of an \textit{adhesive category} \cite{Lack2004}, in which every pushout along a monomorphism satisfies the \textit{van-Kampen square condition}, such that, for every commutative diagram of the form:

\begin{equation}
\begin{tikzcd}
B^{\prime} \arrow[ddd, "f_{h}^{\prime}"] \arrow[dr, "h_B"] & & & A^{\prime} \arrow[lll, "g_h"] \arrow[dl, "h_A"] \arrow[ddd, "f_h"]\\
& B \arrow[d, "f^{\prime}"] & A \arrow[l, "g"] \arrow[d, "f"] &\\
& D & C \arrow[l, "g^{\prime}"] &\\
D^{\prime} \arrow[ur, "h_D"] & & & C^{\prime} \arrow[lll, "g_{h}^{\prime}"] \arrow[ul, "h_C"]
\end{tikzcd},
\end{equation}
such the following sub-diagrams are both pullbacks:

\begin{equation}
\begin{tikzcd}
B^{\prime} \arrow[d, "h_B"] & A^{\prime} \arrow[l, "g_h"] \arrow[d, "h_A"]\\
B & A \arrow[l ,"g"]
\end{tikzcd}, \qquad \text{ and } \qquad
\begin{tikzcd}
A \arrow[d, "f"] & A^{\prime} \arrow[l, "h_A"] \arrow[d, "f_h"]\\
C & C^{\prime} \arrow[l, "h_C"]
\end{tikzcd},
\end{equation}
the pushouts and pullbacks satisfy a compatibility condition. In fact, \textit{selective adhesive} categories \cite{Dixon2013}, \cite{Kissinger2011} allow DPO rewriting to be defined in an even broader class of cases, in which one is working within a full subcategory ${\mathbf{C}^{\prime}}$ of an adhesive category ${\mathbf{C}}$, with embedding functor ${S : \mathbf{C}^{\prime} \to \mathbf{C}}$, with the only condition being that $S$ preserves monomorphisms.

This construction allows us to formalize the notion of \textit{rulial space}, i.e. the space of all possible rewriting rules of a given class (e.g. hypergraph transformation rules, string substitution rules, Turing machine rules, etc.) to be followed between states of a system. For the case of Wolfram model  graph/hypergraph rewriting systems (and other related systems that are known to be special cases of these), this definition can be formalized as:

\begin{definition}
The `rulial space' of Wolfram model systems is  the category of cospans of a DPO rewriting system, defined over the selective adhesive category of Wolfram model hypergraphs.
\end{definition}
The rulial space of Wolfram model systems functorially acquires the structure of a selective adhesive category. Indeed, as a consequence of the concurrency and parallelism theorems of algebraic graph transformation theory  \cite{Ehrig2006}, the rulial space also inherits a natural monoidal structure  \cite{Gorard2020c} (which is, in turn, inherited by all Wolfram model multiway systems, since they are obtained by the composition of certain rules in rulial space).

\subsection{Relation to Other Works }

The first predecessors of the sort of rewriting systems considered in this article   go back to early work on cellular automata and graphical models discussed in \cite{Wolfram2002a}. Those rewriting systems were employed to model a wide range of complex systems across physics, biology and computer science (including patterns of crystal growth, turbulent phenomena, chaotic systems, morphological patterns,  pigmentation of mollusk shells, Turing machines, register machines,  combinators, symbolic substitution systems, tag systems, etc.; just to name a few).  However,  in \cite{Wolfram2020}, it was argued that a   further  generalization  to abstract rewriting systems (i.e., those where rewriting states can be hypergraphs, open graphs, string systems or a formal symbolic language) may be necessary, particularly, for describing  models of physics where rigid notions of spacetime are superseded by geometric structures that themselves are dynamic. To the best of our knowledge, these are among the first proposals suggesting the fundamental role of generic rewriting systems to physics (and our work here is an attempt to formalize such models in the language of  modern homotopy type theory, resulting in a higher categorical framework for investigating pregeometric structures relevant to physics).  

Besides the Wolfram model, another related research program, which makes use of diagrammatic reasoning for describing quantum circuits and quantum processes, is categorical quantum mechanics (CQM) \cite{abramsky2009categorical}, \cite{coecke2018picturing}.  Though the original formulation of CQM is not explicitly based on  rewriting (other recent attempts, besides \cite{Gorard2020c}, \cite{Gorard2021a}, trying to formalize quantum processes within the framework of string rewriting theory can be found in \cite{bonchi2020stringI}, \cite{bonchi2021stringII}), it shares certain similarities with the Wolfram model, in that, it is a framework for expressing diagrammatic process algebras (those based on monoidal categories). This treatment of symmetric monoidal categories as a general language for reasoning about physical systems (with morphisms between objects playing the role of physical transformations between states, and with morphism composition and monoidal composition playing the role of sequential and parallel combination of such processes, respectively) also provides the  mathematical formalization of   physical processes described on Wolfram model multiway systems defined via abstract rewriting rules over arbitrary symbolic expressions.  Some recent   developments to the original Wolfram model can be found in  \cite{Gorard2020c}, \cite{Gorard2021a},  which advanced a formulation in terms of monoidal categories, and which, in turn allowed for a Wolfram model realization of categorical quantum mechanics and a multiway diagrammatic framework for quantum circuits using ZX calculus. In particular, it was shown that Wolfram model multiway systems (whose states are ZX diagrams) serve as  formal embedding spaces of  ZX  processes of CQM, with  multiway rewriting rules precisely corresponding to equational rules of ZX calculus.  

Outside of physics, rewriting systems have had a long history in theoretical computer science.  They have been used  extensively  in algebraic   logic, proof theory and formal programming languages  \cite{Baader1998},  \cite{Bezem2003},  \cite{Dershowitz1990} ,  \cite{Huet1980}.  There exists a large body of relevant work on  graph transformations and graph grammars  \cite{ehrig1999handbook}, \cite{drewes1997hyperedge}, \cite{heckel2006graph}, \cite{heckel2020graph}, including algebraic  graph transformations \cite{corradini1997algebraic}, \cite{hartmut2006fundamentals}. Among the broader classes of rewriting systems that exist in the literature, the ones that interest us here are specific in two respects: (i) we will consider  non-deterministic rewriting (what are called multiway systems), that keep track of all possible rule applications and orderings (consequently, all possible evolution histories);  (ii)  we also retain data of the causal structure of rewriting events (in other words, causal edges in multiway systems are morphisms of a suitable category).  That being said, besides the classical rewriting literature cited above, which has mainly been directed at formal models of computation, let us also point to other recent works in formal rewriting theory.  First, there are the following developments in rule algebras: \cite{sobocinski2020rule}, \cite{behr2021compositionality}, used for stochastic and combinatorial rewriting systems \cite{behr2021rewriting}; and also tracelet analysis  \cite{behr2019tracelets}.  In particular, in recent work \cite{behr2021tracelet}, it was shown that tracelets give rise to decomposition spaces, that is, certain types of simplicial groupoids, and combinatorial Hopf algebras. 
This is relevant when considering concatenation of direct derivations versus parallel execution, as that is only exactly compatible if one considers rules and derivations up to isomorphism, which is what the tracelet Hopf algebra construction  addresses via suitable equivalence relations. In future work, the use of these methods in the Wolfram model will be explored.  
Then, there have also been  developments in the mathematical foundations of  string diagram rewriting, reported in  \cite{bonchi2020stringI},  \cite{bonchi2021stringII}  and   \cite{bonchi2021stringIII}.  Given the wide applicability of the graphical syntax of string diagrams in symmetric monoidal categories, these authors have sought to lay out a thorough mathematical foundation for string diagram rewriting, which unlike term rewriting, poses additional challenges. In particular, they introduce a combinatorial interpretation of string diagram rewriting modulo Frobenius structures, in terms of double-pushout hypergraph rewriting. Since the Wolfram model makes use of both, double-pushout graph and double-pushout hypergraph rewriting, the above works provide the 1-categorical definitions for the kind of rewriting systems we will use here. 
Furthermore, in addition to the above, there is also other relevant work in higher dimensional rewriting theory, and in particular, its incarnation in the theory of polygraphs\footnote{We thank Amar Hadzihasanovic for pointing us to this literature.},  stemming from the pioneering work of Burroni  \cite{Burroni1993}.  Within the framework of polygraph rewriting theory, one can  define  higher rewriting or oriented syzygies that fill confluence  squares in higher dimensions, which form structures of higher categories. These systems can then be investigated with homotopical and  homological tools. Recent reviews on this work can be found in \cite{Guiraud2018},  \cite{Guiraud2019}, and a geometrical perspective based on weak higher categories can be found in \cite{Hadzihasanovic2020}. In contrast, in our work here, we will be concerned with strict higher categories. It would be interesting for future work to pursue the ``categorical weakening" of our results here and relate them to polygraph rewriting.  

Yet another relevant framework for investigating higher structures concerns what are called `hyperstructures'  \cite{baas2016higher}.   These address the general problem of taking local information and properties to extract global information; and have been applied to a variety of complex systems from brains, to languages, to geometries of geometries, to engineering  \cite{baas2019mathematics}, \cite{baas2019philosophy}. In a typical setting, global data is formally extracted from a system's parts by using gluing techniques over presheaves and associated sheaves, or in the categorical settings over Grothendieck topologies and sites. What hyperstructures do is that they extend the notion of Grothendieck topologies, sites and (pre)-sheaves in such a way that gluing is well-defined in a more general context. Presumably, this framework of hyperstructures will be closely related to our constructions of homotopical rewriting systems and will certainly be interesting for future explorations.

As mentioned in the introduction, our work borrows heavily from recent  advances in the foundations of mathematics, particularly homotopy type theory  \cite{Program2013}  and synthetic geometry  \cite{Shulman2016},  \cite{Shulman2017}, \cite{Schreiber2013a}.  Our long-term goals with respect to homotopy type theory and higher categorical structures to physics are:  (i) to explore higher symmetries and spaces, that cannot    readily  be captured by current methods; and (ii)  to seek a constructivist foundation for theories of physics in much the same way that such a foundation is proving fruitful for mathematics itself. From that perspective, this paper is simply a modest attempt at understanding how pregeometric structures and space itself may arise from abstract rewriting systems, and more generally, from the combinatorics of discrete computational building blocks. For future work, it may be informative to also express particles, their interactions and symmetries in this type-theoretic framework. A related but different approach has been pursued in  \cite{Schreiber2012},  \cite{Schreiber2013a},  where quantum field theories with  higher gauge symmetries have been studied using homotopical methods. Likewise, homotopical pre-quantum geometries have been discussed in  \cite{schreiber2016higher}  (see also  \cite{isham2000some} for an early predecessor of such methods).  All of these works differ from the pregeometric homotopical structures we consider in this article. Of course, the main difference with these  works   is that they seek to extend the current framework of quantum theory and  symmetries to higher symmetries (and consequently to Lie algebroids), while our work here is an attempt to seek a foundational understanding of existing physical structures from discrete underpinnings. On the other hand, homotopical methods naturally lead one to the unifying formalism of topos theory and constructivist mathematics, and that is the conceptual similarity with our work here.


\section{Multiway Rewriting Systems Formalized as a Type Theory}

In algebraic and combinatory logic, rewriting systems are commonly used for theorem proving, where a sequence of rewriting terms corresponds to a proof of  a proposition. In fact, these systems are formal models of computation, equivalent to Turing machines.  In practice, type theory can be used as an extremely powerful meta-language to express various formal computational systems, including rewriting systems; and that is what we will use here. 

The archetypical constructions of the Wolfram model are multiway rewriting  systems:  non-deterministic rewriting systems equipped with a causal structure.   Compared to pure rewriting systems, these additionally preserve causal information as chains of partial order in the multiway graph. Moreover,  multiway systems admit all permissible applications of  rewriting rules in all possible orders, leading to parallel causally-ordered branches of rewriting sequences. Furthermore, Wolfram model multiway systems may be realized over graphs, hypergraphs or strings. In this section, we show how these systems can be expressed in type theoretic syntax.  This general approach of seeking a constructivist formalization of  fundamental structures has strong parallels in recent work at the foundations of  mathematics and computer science.

\subsection{Type Constructors for Multiway Rewriting Systems} 

A multiway system in the Wolfram model is a network of states connected by causal edges determined by  rewriting rules that inductively generate the complete (either terminating or infinite) multiway system.  Here, we will show that Wolfram model rewriting rules can be formally expressed as computational rules of type theory. This representation of multiway systems enables one to inductively construct abstract mathematical structures starting from Wolfram model multiway systems in much the same way that one does within the formalism of type theory. Following that, we make the case that multiway systems themselves can be thought of as types; multiway states, as terms of a type; and multiway rules, as computational type constructors. We first show the type construction for a string rewriting multiway system and subsequently for one involving hypergraph rewriting. 

\subsubsection{String Rewriting Multiway Systems}

\begin{figure}[h]
\centering
\includegraphics[width=0.9\textwidth]{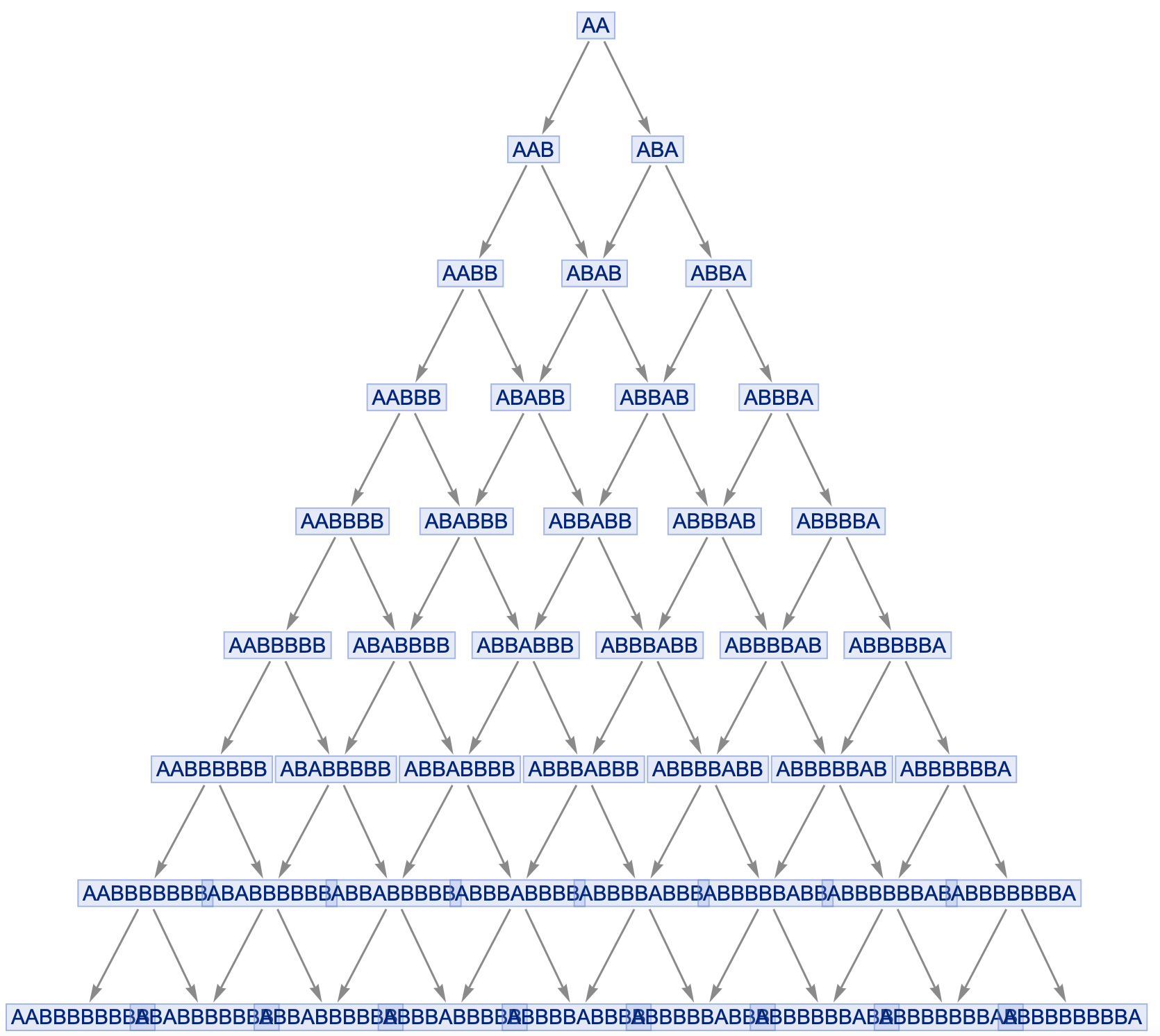}
\caption{The multiway evolution graph corresponding to the first 8 steps in the non-deterministic evolution history of the string substitution rule ${A \to AB}$, starting from the two-character initial condition $AA$.}
\label{f3.1.1}
\end{figure}

Let us begin with a prototypical example of a string rewriting rule such as ``A" $\to$ ``AB", which generates the multiway network shown in Fig. \ref{f3.1.1}. Even though this is a specific example, what we illustrate here is generic to all multiway systems.  In examples such as this, one is given an initial state. Repeated application of the rule(s) in all possible ways to that state and all subsequent states will generate the entire multiway system, and causal chains of terms thus generated will correspond to proofs (or proofs of equality in type theory) that a given term implies another.

In practice, type theory can be used as an extremely powerful meta-language to formalize  models of computation, including rewriting systems. 
In type theory, constructors constitute a package of rules defining types. There are four basic classes of rules in type theory: (i) `Formation' rules such as  $\Gamma \vdash A \,\, \mathsf{Type} $  declare the existence of a type. These are also called type judgements; (ii) Then there are `Introduction' rules as  $\Gamma \vdash  a : A $. These are called term judgements, indicating terms associated to a type. In the Wolfram model, when expressed as a type theory, type judgements will correspond to declaration of a multiway system and term judgements will refer to declaring states (nodes) of a multiway graph (which is explicitly shown below); (iii)   `Elimination' rules in type theory declare how to use and /or substitute newly introduced terms. These will correspond  to the  familiar relational rules of the Wolfram model; and (iv) Finally, there are  `Computational' rules which declare equivalences between terms of a type. These will refer to state equivalences in the Wolfram model.  Additionally, for what follows, we will be working with dependent type theory, where types themselves may depend on other types.
 
Below we state the explicit type constructors for the Wolfram model multiway system in Fig. \ref{f3.1.1}.  This is the explicit algorithmic description of multiway systems that any automated proof system could reason with. 

To begin, we will first need to construct a monoid $S$ using the following formation, introduction, elimination and computation rules:
\begin{eqnarray}
  \inferrule{\Gamma \vdash S \,\, \mathsf{Type} }{\Gamma  \vdash \mathbf{1} : S \quad  \Gamma  \vdash s_1 : S \quad \cdots \quad \Gamma  \vdash s_n : S}  
\end{eqnarray}
\begin{eqnarray} 
 \inferrule{\Gamma  \vdash x : S \quad  \Gamma  \vdash y : S}{\Gamma  \vdash x \, y : S}
 \end{eqnarray}
\begin{eqnarray} 
 \inferrule{\Gamma  \vdash x : S \quad  \Gamma  \vdash \mathbf{1} : S}{\Gamma  \vdash x \, \mathbf{1} \equiv x \quad \Gamma  \vdash \mathbf{1} \, x \equiv x}
 \end{eqnarray}
where $\mathbf{1}$ is the identity element and the $s_i$ are generators. This package of rules constructively realizes a monoid type, and hence, is also referred to as a type constructor for monoids. The product is given by concatenation. Elements of $S$ are strings built out of a collection of symbols $s_i$. In what follows, the terms of our (yet to be constructed) multiway system will be generated from the monoid $S$.

Next, we need to specify rewriting maps $\{ {\cal R} : S \to S  \}$,  that constitute rules for transforming an initial string to a final string $\{ {\cal R} (s_{ini}) : s_{ini} \to s_{fin}  \}$. These are simply function types and can be expressed using the standard type constructors for functions. The maps $\{ {\cal R}  \}$ are, in fact, Wolfram model rules. 

We also need to construct a dependent type ${\cal M}[S]$ whose terms $m[s]$ constitute states of the multiway system, and an inclusion map $ i : {\cal M}  \hookrightarrow S $ from ${\cal M}$ to $S$.  All of that can be specified via the following type constructions:


\begin{eqnarray} 
  \inferrule{\Gamma, s : S \vdash {\cal M} \,\, \mathsf{Type} }{\Gamma, s : S   \vdash m [s] : {\cal M} }   
\end{eqnarray}
\begin{eqnarray}   
  \Gamma \vdash ({\cal M} \hookrightarrow S) \,\, \mathsf{Type}  
\end{eqnarray}
\begin{eqnarray}   
  \inferrule{\Gamma  \vdash  m : {\cal M}  }{\Gamma    \vdash \lambda m.i : {\cal M} \hookrightarrow S }   
\end{eqnarray}
\begin{eqnarray}   
  \inferrule{\Gamma  \vdash i : {\cal M} \hookrightarrow S  \quad  m : {\cal M} }{\Gamma    \vdash i (m) : S }   
\end{eqnarray} 

Furthermore, we need the inverse function of $i$, that gives a term in ${\cal M}$ for a given term $s$ in $S$. This can be defined via a function type, which  satisfies the following computational rule:

\begin{eqnarray} 
  \inferrule{\Gamma, m : {\cal M}, s : S  \vdash i (m) \equiv s }{\Gamma  \vdash m [s] \equiv i^{-1} (s) }   
\end{eqnarray}  

We now have an algorithmic procedure for constructing terms in ${\cal M}$ from terms in the monoid $S$. In general, a term $m$ can depend on $s$ in many number of ways. Here, we further specify that $m [s] \equiv s$. This will yield multiway states as a collection of strings from a symbol set. Note that, in  type theory, ${\cal M}$ and $S$ are necessarily distinct types, with the former  being  dependent on the latter. Hence, it was necessary to go through this explicit procedure of first constructing $S$ as a bag of strings, from which one can then pick up those that yield states of the above multiway system.

Given the above, the following introduction rule, declares the initial states of the multiway system:
\begin{eqnarray}
  \Gamma \vdash  m_{in} :  {\cal M}   
\end{eqnarray}

Note that any term $\Gamma \vdash  s : S$ is, by construction,  a disjoint union of the generators (and identity) of $S$ 
\begin{eqnarray}
 s \equiv  \coprod_{l = 1}^{|s|} s_{k_l} 
\end{eqnarray}
where $|s|$ denotes the cardinality of the word $s$. For what follows, it will be convenient to pad any $s$ with the identity element as so:  $  \mathbf{1} \cdot s \cdot \mathbf{1}  $.  

We now need to set-up type theoretic rules for instantiating string rewriting. Recall the rewriting maps we constructed above. We now need to specify an algorithmic type theoretic procedure for explicitly realizing this rewriting (in a way that an automaton could carry out these operations). 

We first need a map that expresses the padded term $  \mathbf{1} \cdot s \cdot \mathbf{1}  $  in terms of all of its tri-partitions: 
\begin{eqnarray}
 \Gamma \vdash  {\cal P}_3 : \mathbf{1} \cdot S \cdot  \mathbf{1}  \to \{ S \coprod S \coprod S \}  
\end{eqnarray} 
\begin{eqnarray}
 \Gamma \vdash  {\cal P}_3 \left( \mathbf{1} \cdot s \cdot \mathbf{1} \right) \,\,  \equiv  \,\, \{ s_a \coprod s_{b} \coprod s_c \}  
\end{eqnarray}
This map exists due to the fact that any $s$ is, by construction,  a disjoint union of generators (identity elements) of $S$. The partitioning operation simply entails all possible 3-groupings. We can compactly express this via the map $\widetilde{{\cal P}_3}$, where the padding and partitioning are all done at once:
\begin{eqnarray}
 \Gamma \vdash \widetilde{{\cal P}_3} : S \to \{ S \coprod S \coprod S \}  
\end{eqnarray} 
\begin{eqnarray} 
 \Gamma \vdash  \widetilde{{\cal P}_3} (s) : \{ S \coprod S \coprod S \}  
\end{eqnarray} 
\begin{eqnarray} 
 \Gamma \vdash  \widetilde{{\cal P}_3} (s) \equiv \{ s_a \, s_b \, s_c \}  
\end{eqnarray}

Lifting these maps to ${\cal M}$, entail the following rules:
\begin{eqnarray}
  \Gamma \vdash  \widehat{{\cal P}_3} :  {\cal M}  \to \{ {\cal M} \}  
\end{eqnarray} 
\begin{eqnarray}
 \Gamma \vdash  \widehat{{\cal P}_3} ( m [s] ) \equiv  \{ m [ s_a \, s_b \, s_c ] \}  
\end{eqnarray} 

With these ingredients in place, term rewriting over states of  ${\cal M}$ can be instantiated via the following type theoretic rule:
\begin{eqnarray}
  \inferrule{\Gamma  \vdash  s : S, \, \, m [s] : {\cal M}, \, \, s_{ini} : S, \, \, s_{fin} : S, \, \,  \widetilde{{\cal P}_3} (s) \equiv \{ s_a \, s_b \, s_c \}, \, \, \widehat{{\cal P}_3} ( m [s] ) \equiv  \{ m [ s_a \, s_b \, s_c ] \}    }{\Gamma  \vdash  \mathsf{Conditional}\!\left\{ \left( s_b \equiv s_{ini}, \, \, m [ s_a \, ( s_{fin} / s_b ) \, s_c ], \, \,  m [ s_a \, s_b \, s_c ]  \right) \right\}  }  
\end{eqnarray} 
where the type '$\mathsf{Conditional}$' operates as an 'If...then...else' statement, read as, 'If $s_b \equiv s_{ini}$, then $s_b$ is replaced by $s_{fin}$, else left unchanged. The braces $\{...\}$ denote that this operation is applied to every tri-partition of $s$ such that all possible rewritings of a given state $m[s]$ in ${\cal M}$ are instantiated. 

Besides terms $m[s]$ of the type ${\cal M}$, corresponding to the collection states, our multiway system also contains additional data: that associated to  causal ordering of rewriting events (or partial ordering of terms). This data is captured via maps $\widehat{{\cal R}}$, tied to cases from the conditional rule above whenever a substitution of $s_b$ took place. This can be expressed as:
\begin{eqnarray}
\Gamma \vdash  \widehat{{\cal R}}_{s_a \, s_b \, s_c} : {\cal M} \to  {\cal M}  
\end{eqnarray} 
\begin{eqnarray}
\Gamma \vdash   \widehat{{\cal R}}_{s_a \, s_b \, s_c} : m [ s_a \, s_b \, s_c ]\text{\Large{$|$}}_{s_b \equiv s_{ini}}  \to \, \, m [ s_a \,  s_{fin}  \, s_c ]  
\end{eqnarray}

Finally, we need one more rule to implement state equivalences in the multiway system. This is necessary when one has two distinct rewriting events which result in equivalent final state (i.e., corresponding to identical words of the monoid $S$). These equivalent final states may well be initiated by two distinct initial states or even the same initial states but with two different rewriting rules. This identification of equivalent states is what results in mergers between causal edges in the multiway graph in  Fig. \ref{f3.1.1}. rewriting events themselves lead to splits in the graph, state equivalences lead to mergers. The type theoretic rule to implement  state equivalence is
\begin{eqnarray}
  \inferrule{\Gamma  \vdash  s_i : S, \, \, s_j : S,  \, \,  m [s_i] : {\cal M},  \, \,  m [s_j] : {\cal M}  }{\Gamma  \vdash  \mathsf{Conditional}\!\left\{ \left( s_i \equiv s_j, \, \, m [s_i]  \equiv m [s_j],   \, \,  ---  \right) \right\}  }  
\end{eqnarray} 
where the "$---$" denotes a null operation.

This concludes the full set of type-theoretic rules or type constructors needed to fully specify a string substitution multiway system $\left( {\cal M}, \widehat{{\cal R}}  \right)$ as a rewriting system with causal structure. The multiway graph in  Fig. \ref{f3.1.1}  provides a presentation of a multiway type with all its terms.

\subsubsection{(Hyper)Graph Rewriting Multiway Systems}
 
The corresponding type constructors for Wolfram model multiway systems involving graph and  hypergraph rewriting are effectively expressed as extensions to those described above for string rewriting. In particular, these are treated as set substitution systems, replacing string substitutions.  The main issues we then have to specify is how the states of this multiway system are defined and what state equivalences mean in this context. Apart from that, the rest of the algorithm will be similar to that of string substitution multiway systems described at length above. 

First, let us recall the set-theoretic representation of hypergraphs used in the Wolfram model:  Given a vertex set $V$ with $n \in {\mathbb N}$ elements $v_i$ and $1 \leq i \leq n$, we will simply use vertex labels $i$ to denote vertices or nodes themselves. Then a directed edge between from vertex  $i$ to $j$ is the ordered list  $\{ i, j \}$.  Likewise, a hyperedge consisting of $m$ nodes (which we will call the `arity'  of the edge) can be expressed as a finite collection or list of ordered relations  $\{ i_1, \cdots, i_m \}$ with $i_a \in {\mathbb N}$, as shown in Figure \ref{fig:Figure1}.     
That being said, a Wolfram model hypergraph can be expressed as a disjoint union of ordered lists of natural numbers; for example,  ${\left\lbrace \left\lbrace 1, 2, 3 \right\rbrace, \left\lbrace 3, 4, 5 \right\rbrace \right\rbrace}$ (see Figure \ref{fig:Figure1}).  In other words, given hyperedges $E_l$, we write a coproduct of ordered lists 
\begin{equation}
\coprod_{l = 1}^k  E_l 
\end{equation}
to denote a hypergraph consisting of $k$ hyperedges defined from a vertex set   $V =  \{ 1, \cdots, n \}$.  These constitute the states of a hypergraph multiway system.  
 
How are these multiway states expressed type theoretically? The relevant type constructors for this are the dependent coproduct type and a monoid type. For the former, one can simply invoke the usual  coproduct constructor (this has been nicely detailed in  \cite{Shulman2017}); whereas, the latter, would be similar to  the monoids we constructed for the string substitution multiway above. Each $E_l$  is a monoid over a finite generating subset of ${\mathbb N}$ with the identity element simply being the null element $\phi$.  With these type constructors defining states of the hypergraph multiway system, most of the algorithmic procedure for constructing string-based multiway system (including tri-partitioning the states and checking for matching lists) follows straightforwardly.  The rewriting rules for this multiway system are now  hypergraph transformations rules as exemplified in  Figure \ref{fig:Figure2} (in fact, these rules are a family of rules that apply for several subsets of naturals that match the pattern).  The evolution  history of  this system is shown in  Figures \ref{fig:Figure3} and  \ref{fig:Figure4}, and the multiway graph in Figure \ref{fig:Figure5}.  The other issue concerns state equivalences. In the Wolfram model these are specified via approximate hypergraph isomorphisms, which are based on extensions of standard error-tolerant graph matching algorithms \cite{bunke1998error},   \cite{gorard2016uniqueness}. The type theoretic representation of these algorithms are straightforward but tedious. 

Even though, for practical purposes it has been convenient to program hypergraph rewriting as set substitution systems, it is worth pointing out that more formally hypergraph rewriting can in fact be expressed as a special case of string diagram rewriting (not to be confused with character string rewriting discussed above). We point the reader to recent advances in this direction discussed in \cite{bonchi2020stringI}.

\subsection{Groups, Categories and Groupoids from Rewriting Systems} 
 
In type theory, the structures one constructs are merely syntactic entities. It is only upon being interpreted within a suitable classifying category that these constructs express familiar mathematical objects and their properties. The semantics of types correspond to universal properties internalized within an appropriate category. Likewise, multiway rewriting rules and multiway graphs  are also syntactic structures, which can be interpreted within an appropriate category or groupoid. In type theory, constructors correspond to rule packages associated to mathematical operations or objects such as product types, function types, natural number types, etc. A well-ordered application of these constructors corresponds to respective mathematical objects (interpreted within a category). For example, with type constructors corresponding to a free product, a generator, an identity and an inverse (along with rules for associativity) one can generate the syntax of a free group. The type-theoretic rules for this are expressed as follows:

\begin{eqnarray} 
\inferrule{\Gamma \vdash G \,\, \mathsf{Type} }{\Gamma  \vdash e : G \quad  \Gamma  \vdash g_1 : G \quad \cdots \quad \Gamma  \vdash g_n : G}  
\end{eqnarray} 
\begin{eqnarray}
 \inferrule{\Gamma  \vdash x : G \quad  \Gamma  \vdash y : G}{\Gamma  \vdash g [ x, \, y ] : G}
\end{eqnarray} 
\begin{eqnarray}
 \inferrule{\Gamma  \vdash g [ x_1, \, y_1 ] : G \quad  \Gamma  \vdash g [ x_2, \, y_2 ] : G}{\Gamma  \vdash g [ g [ x_1, \, y_1 ], \, g [ x_2, \, y_2 ] ] : G}
\end{eqnarray} 
\begin{eqnarray}
 \inferrule{\Gamma  \vdash x : G \quad  \Gamma  \vdash e : G}{\Gamma  \vdash g [ x, \, e ] \equiv x \quad \Gamma  \vdash g [ e, \, x ] \equiv x}
\end{eqnarray} 
\begin{eqnarray}
 \inferrule{\Gamma  \vdash x : G \quad  \Gamma  \vdash x^{-1} : G}{\Gamma  \vdash g [ x, \, x^{-1} ] \equiv e \quad \Gamma  \vdash g [ x^{-1}, \, x ] \equiv e}
\end{eqnarray} 
\begin{eqnarray}
 \inferrule{\Gamma  \vdash x : G \quad  \Gamma  \vdash y : G \quad  \Gamma  \vdash z : G}{\Gamma  \vdash g[x, \, g[y, \, z] ] \equiv g[ g[x, \, y], \, z]  } 
\end{eqnarray} 

On the other hand, the corresponding multiway graph associated to a free group is shown in Fig. \ref{f3.2.1}. In this case, the axioms of the group have been defined as term  rewriting rules:
\begin{figure}[ht]
\centering
\includegraphics[width=0.6\textwidth]{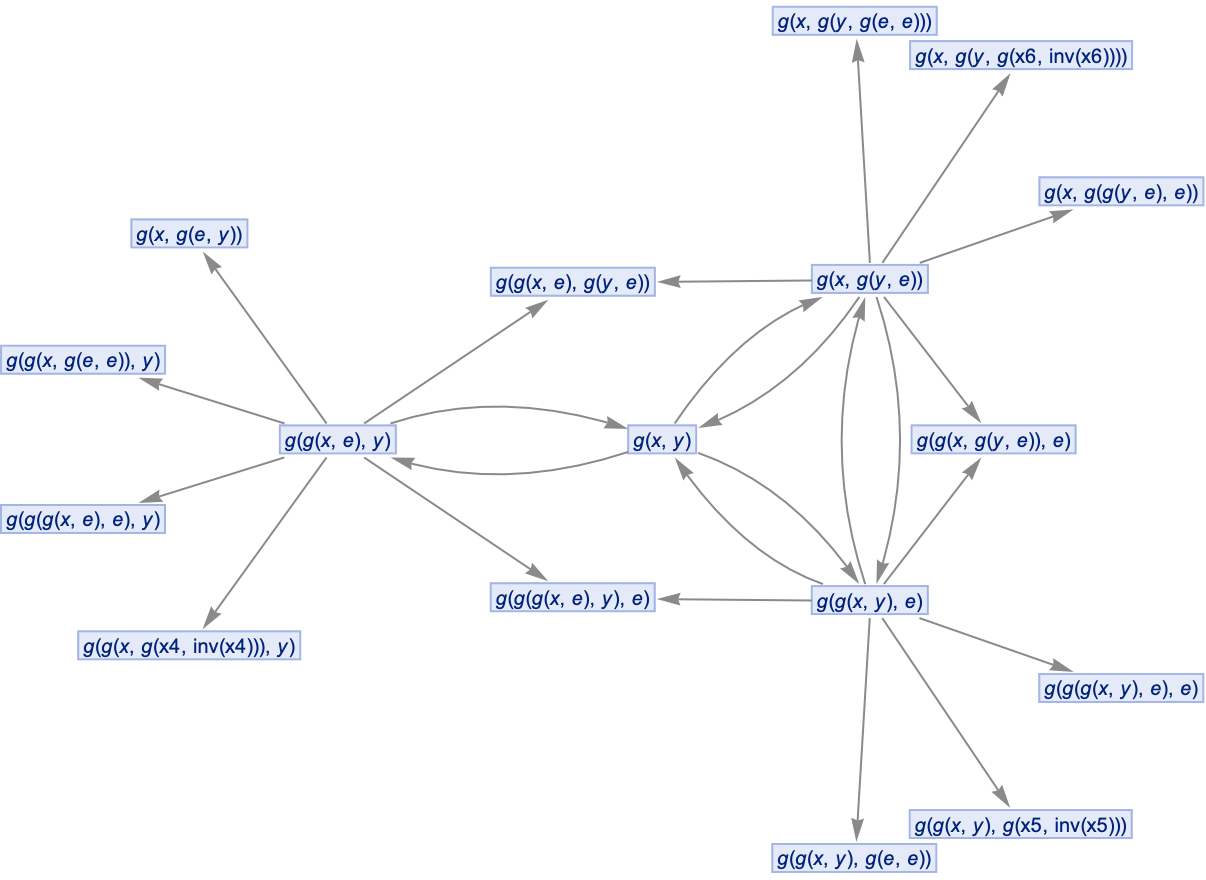}
\caption{The multiway states graph (i.e. a variant of a multiway evolution graph in which cycles are permitted) corresponding to the evolution of the multiway operator system for the axioms of group theory, defined by the term rewriting rules ${\left\lbrace g[x, g[y, z]] \to g[g[x, y], z], \,\,\,\,  g[g[x, y], z] \to g[x, g[y, z]] \right\rbrace}$ (associativity), ${\left\lbrace g[a, e] \to a, \,\, a \to g[a, e] \right\rbrace}$ (right identity rules, plus corresponding rules needed for the left identity) and ${\left\lbrace g[a, \mbox{inv}[a]] \to e, \,\, e \to g[a, \mbox{inv}[a]] \right\rbrace}$ (rules for the right inverse, plus corresponding rules needed for the left inverse).}
\label{f3.2.1}
\end{figure}

Analogous to type-theoretic constructions, the multiway graph of the Wolfram model is capable of diagrammatic presentations of various mathematical structures. As concrete examples of multiway graphs being used for presenting   mathematical structures, let us consider two more constructions which will be relevant for what follows: that of (i) a small category and (ii) a groupoid. 

\begin{figure}[ht]
\centering
\includegraphics[width=0.32\textwidth]{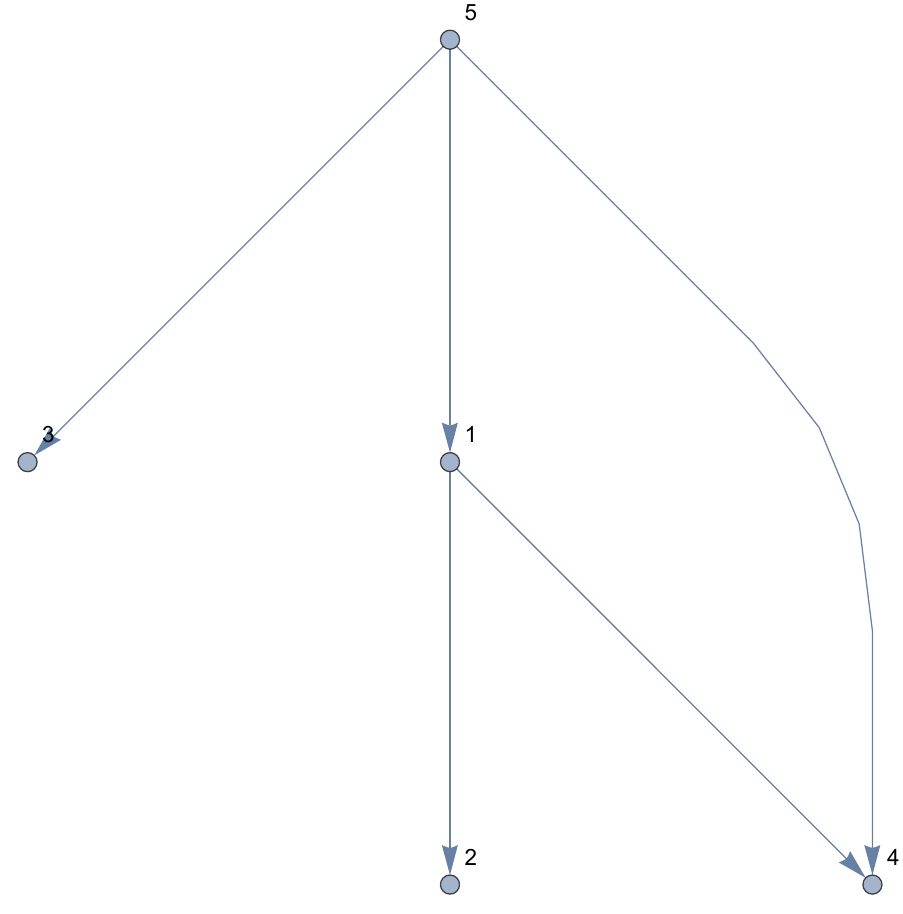}\hspace{0.1\textwidth}
\includegraphics[width=0.5\textwidth]{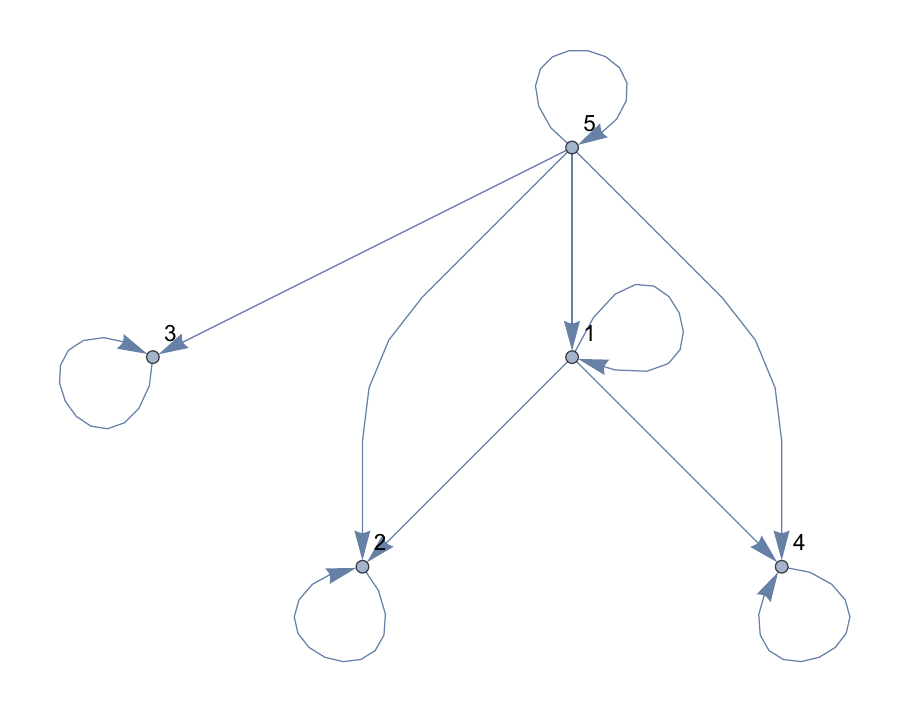}
\caption{An initial graph representing objects and morphisms (shown on left-hand side) and it's categorified version (right-hand side) obtained by graph rewriting using Wolfram model rules.}
\label{f3.2.2}
\end{figure}

A small category is by definition a finite collection of objects (which can be represented as nodes), along with a collection of morphisms between those objects (representable as directed edges), and with the property that   morphisms are associative and reflexive (where self-loops will denote identity morphisms for each object). Diagrammatically, this is simply a directed graph with transitive closure and self-loops. For example, let us consider the  directed graph $G$ shown on the left-hand side of Fig. \ref{f3.2.2}. One can impose associativity by computing the transitive closure of this graph. Following that, one can impose reflexivity by adding self-loops at nodes. The result is shown on the right-hand side of Fig. \ref{f3.2.2}.  

The transformed object on the right-hand side can be  interpreted as a small category (``small" due to a finite set of objects). Moreover, the operations stated above can be explicitly expressed as the following Wolfram model rules applied to any initial graph $G$: 
\begin{eqnarray}
\left\lbrace \left\lbrace a, b \right\rbrace, \left\lbrace b, c \right\rbrace \right\rbrace &\to& \left\lbrace \left\lbrace a, b \right\rbrace, \left\lbrace b, c \right\rbrace, \left\lbrace a, c \right\rbrace \right\rbrace  \\
\left\lbrace \left\lbrace a \right\rbrace \right\rbrace &\to& \left\lbrace \left\lbrace a, a \right\rbrace  \right\rbrace 
\end{eqnarray}
where the first rule refers to transitivity and the second to reflexivity.

\begin{figure}[ht]
\centering
\includegraphics[width=0.3\textwidth]{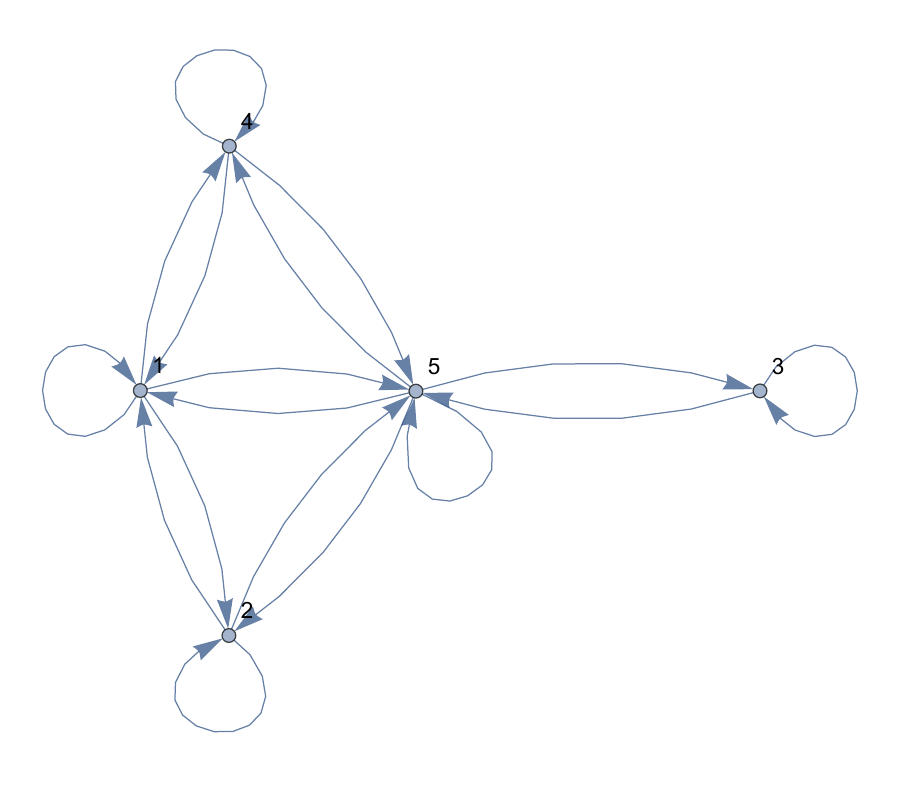} \hspace{0.1\textwidth} 
\includegraphics[width=0.55\textwidth]{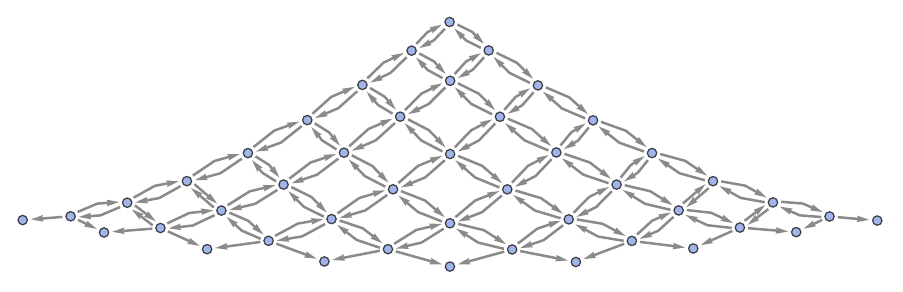}
\caption{Examples of groupoids constructed from graph rewriting using Wolfram model rules. The l.h.s. shows the groupoid resulting from the small category shown in Fig. \ref{f3.2.2}. The r.h.s. shows a multiway system as a groupoid after including invertible Wolfram model rules (besides transitivity and reflexivity, which are not explicitly shown here for clarity of presentation).}
\label{f3.2.3}
\end{figure}

Similarly, a groupoid being a special case of a category in which all morphisms are invertible (and hence isomorphisms), can be  diagrammatically represented  by adding two-way directed edges, indicating isomorphisms, to the graph of the small category. This is shown on the left-hand side of Fig. \ref{f3.2.3}. 

The additional Wolfram model rule ensuring isomorphisms is as follows:
\begin{eqnarray}
\left\lbrace \left\lbrace a, b \right\rbrace \right\rbrace &\to& \left\lbrace \left\lbrace b, a \right\rbrace  \right\rbrace 
\end{eqnarray}
 
Hence, starting from any initial graph (or more generally hypergraphs) $G$ and rewriting with appropriate Wolfram model rules, one can then generate multiway systems corresponding to various mathematical objects such as groups, categories and groupoids. The fact that many algebraic objects have elegant diagrammatic representations is well known. Our emphasis in this section was about demonstrating how these constructions can be realized as multiway rewriting systems that are generated from simple rewriting rules.

 \subsection{A Comment on Homotopies as Proofs of Equality in Multiway Systems}

Homotopy Type Theory (HoTT) is an augmentation of type theory (more specifically, of Per Martin-L\"of's intuitionistic type theory) with one key additional axiom, namely Vladimir Voevodsky's axiom of univalence    \cite{Program2013},  \cite{Ahrens2021}.     The key philosophical idea underpinning homotopy type theory is an extension of the `propositions as types' interpretation of the Curry-Howard correspondence, in which types  correspond to topological spaces (more generally, homotopy spaces), terms of a given type correspond to points in those spaces, proofs of equality between terms of a given type correspond to paths connecting points, proofs of equality between proofs corresponding to higher homotopies between associated paths, etc. In other words, homotopy type theory is a way of endowing type theory with a kind of inbuilt homotopy  structure, where constructing proofs of equality (between proofs) is tantamount to constructing associated homotopy complexes within types.  

In HoTT, one can in fact consider the infinite hierarchy containing all possible higher-order homotopies, corresponding to an infinite hierarchy of all possible higher-order proofs of equality. Types equipped with this structure are formally identified as $\infty$-groupoids. Having formalized multiway rewriting systems  as types, we now proceed to show how the full machinery of higher homotopies and $\infty$-groupoids can be realized by Wolfram model multiway systems.

\clearpage 

\section{Higher Homotopies in Multiway Rewriting Systems as $n$-Fold Categories}

Here we demonstrate how higher homotopies can be constructed using multiway rewriting systems, and introduce a systematic algorithm for identifying rewriting rules that give rise to these homotopies. Subsequently, we prove that a multiway system equipped with homotopies up to order $n$ may be formalized as an $n$-fold category, such that the infinite limit of this higher-order multiway system yields an ${\infty}$-groupoid (upon the admission of invertible rewriting rules). This section extends beyond initial work reported in  \cite{arsiwalla2021homotopy},  especially  parts concerning $\left(\infty, 1\right)$-categories. 

\subsection{An Algorithmic Framework for Higher Homotopies }

Once again, let us consider rewriting systems on character strings. For instance, specified by the rule ${A \to AB}$, generating the multiway evolution graph shown in Fig. \ref{f3.1.1}. We have deliberately selected this minimal example, although the basic algorithm described here applies in principle to  all classes of multiway systems (an analogous algorithmic procedure works for hypergraph or other multiway systems). Every path in this multiway evolution graph corresponds to a proof (in HoTT, paths are proofs of equality between terms).  Shown on the left of Figure \ref{fig:Figure3}  is a  particular case of a proof of  ${AA \to ABBBABBB}$, subject to the axiom ${A \to AB}$.   A key feature of homotopy type theory is that one can consider multiple paths connecting the same pair of vertices, which may be interpreted type-theoretically as corresponding to the existence of multiple proofs of equality, as shown on the right of Figure \ref{fig:Figure3}.  
Additionally, in HoTT one can conceptualize the notion of proofs of equivalences between proofs in terms of homotopies between paths  \cite{Program2013}. In our multiway systems, these will correspond to paths between paths.

\begin{figure}[ht]
\centering
\includegraphics[width=0.4\textwidth]{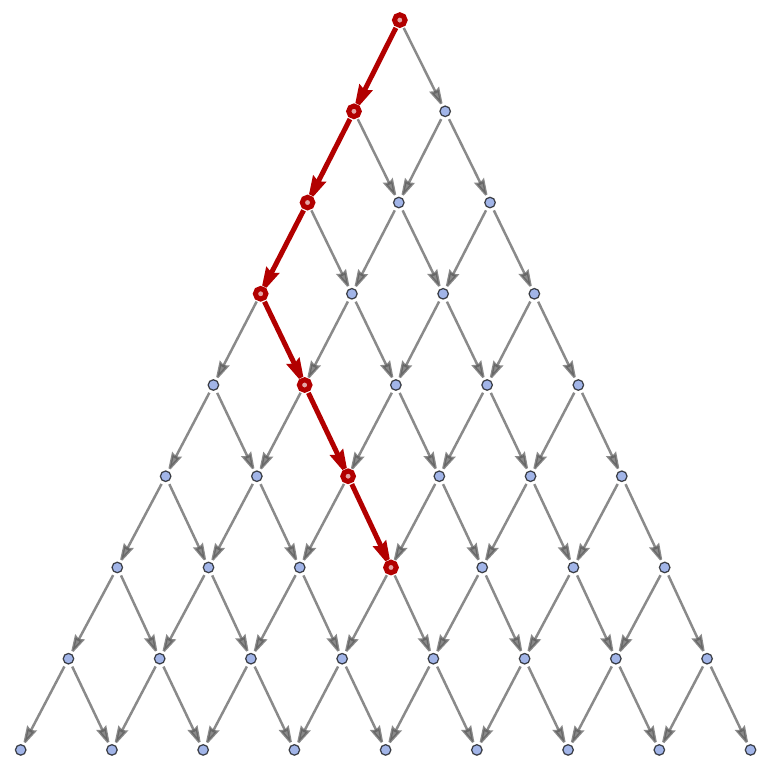}\hspace{0.1\textwidth}
\includegraphics[width=0.4\textwidth]{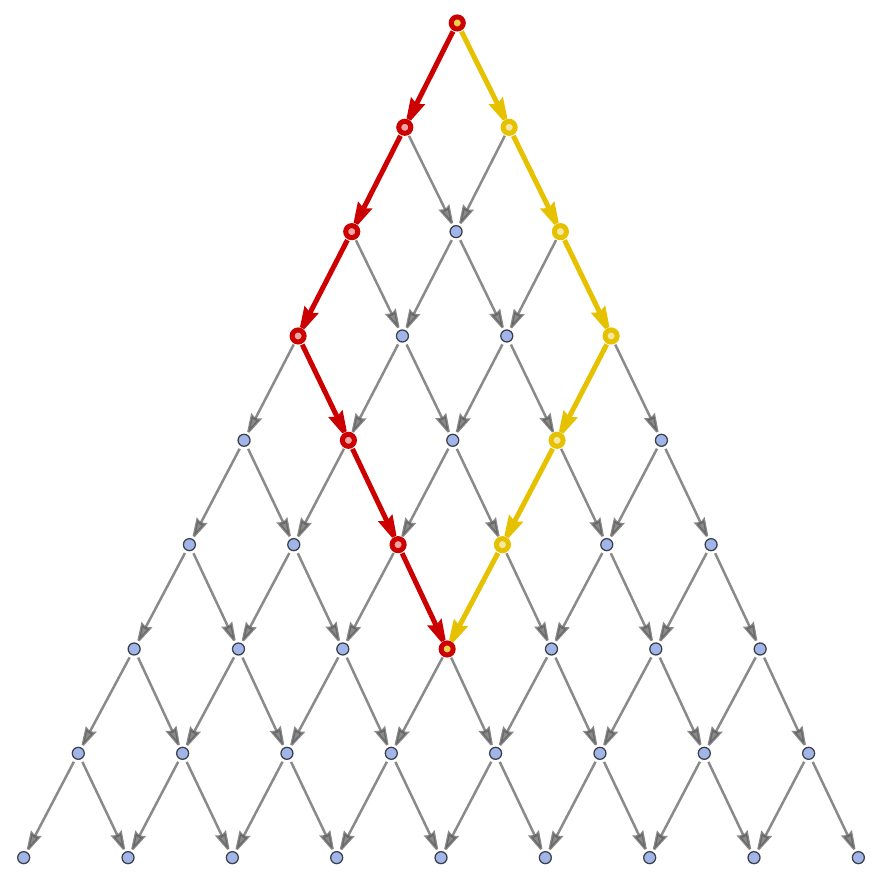}
\caption{On the left, a multiway evolution graph with a red highlighted path between vertices ${AA}$ and ${ABBBABBB}$, corresponding to a proof of  ${AA \to ABBBABBB}$, subject to the string substitution axiom ${A \to AB}$. On the right, the same multiway evolution graph showing multiple highlighted paths (red and yellow) between vertices ${AA}$ and ${ABBBABBB}$, illustrating the existence of multiple proofs of  ${AA \to ABBBABBB}$.}
\label{fig:Figure3}
\end{figure}

In the discrete setup of a multiway evolution graph, a concrete realization of the homotopy map between paths can be given by simply introducing additional  edges mapping vertices from one path to corresponding vertices on the other path,  as shown in Figure \ref{fig:Figure4}. However, this approach to introducing homotopy maps is somewhat ad hoc; it is far more natural to introduce these new mappings between multiway vertices as a new set of rewriting rules that can be appended to the original multiway system rules (shown in Figure \ref{fig:Figure5}). For the particular system at hand, these additional rewriting rules are:

\begin{multline}
\hspace{2cm}  \left\lbrace AAB \to ABA, \,  AABB \to ABBA, \, AABBB \to ABBBA, \right.\\  
\left. ABABBB \to ABBBAB, \,  ABBABBB \to ABBBABB \right\rbrace  \hspace{3cm}
\end{multline}
which yields the system shown on the right-hand side of Figure \ref{fig:Figure5}. 


\begin{figure}[ht]
\centering
\includegraphics[width=0.5\textwidth]{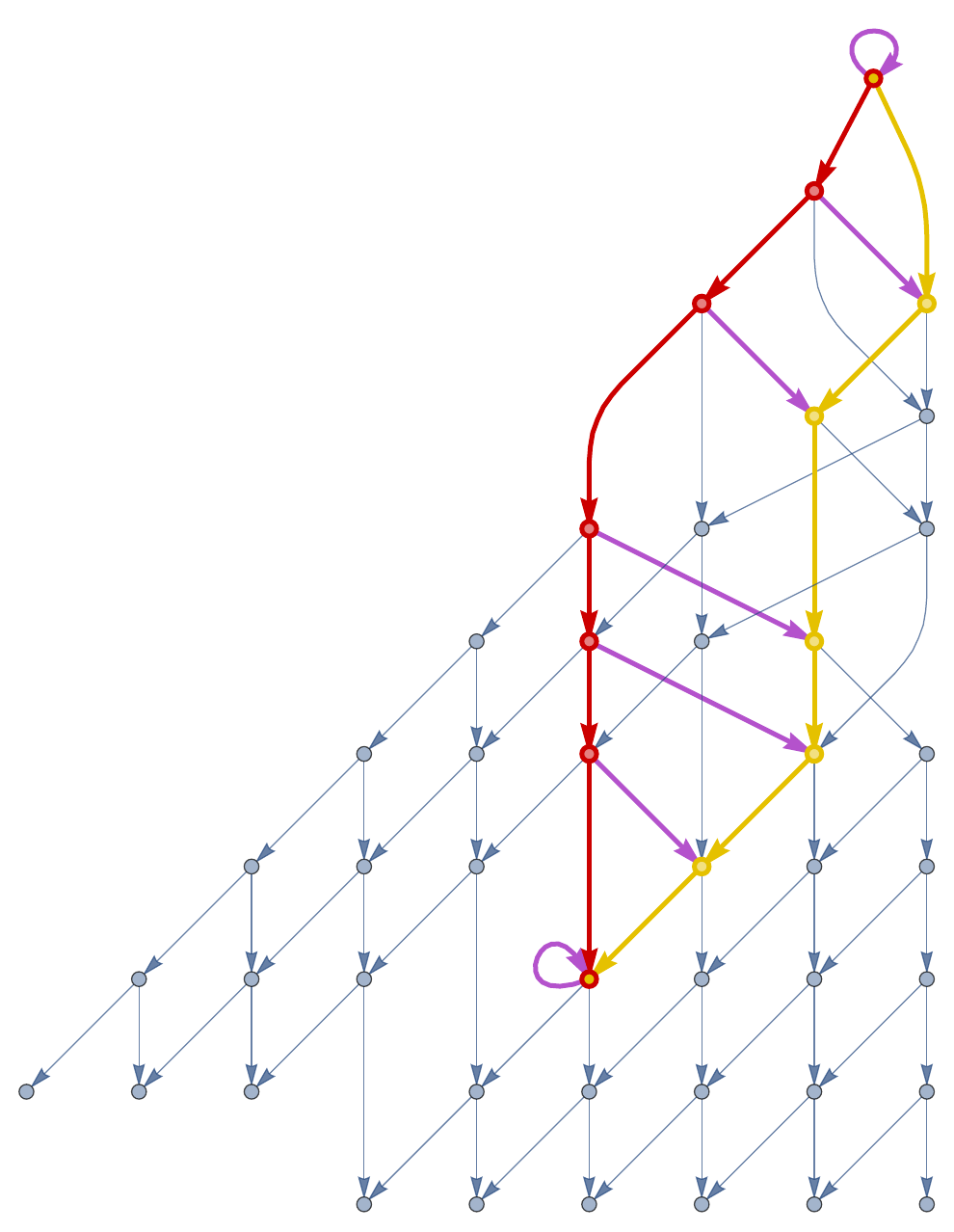}
\caption{A multiway evolution graph with purple highlighted paths between from vertices on the red path to corresponding vertices on the yellow path, interpreted as homotopy maps between the two associated proofs of  ${AA \to ABBBABBB}$.}
\label{fig:Figure4}
\end{figure}

\begin{figure}[ht]
\centering
\includegraphics[width=0.4\textwidth]{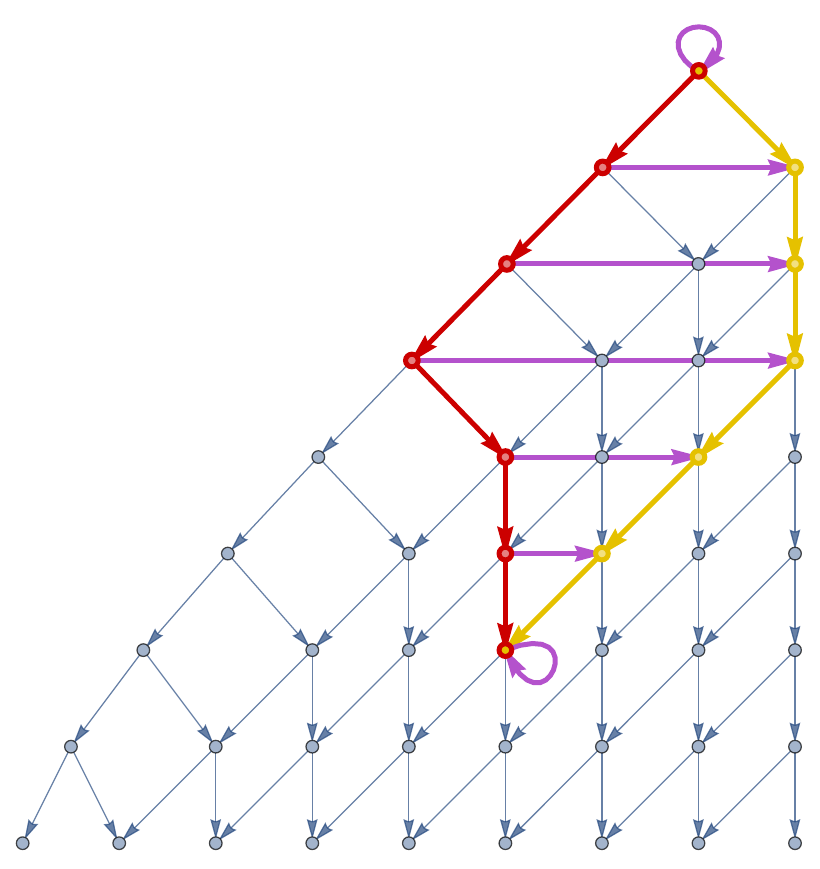}\hspace{0.1\textwidth}
\includegraphics[width=0.4\textwidth]{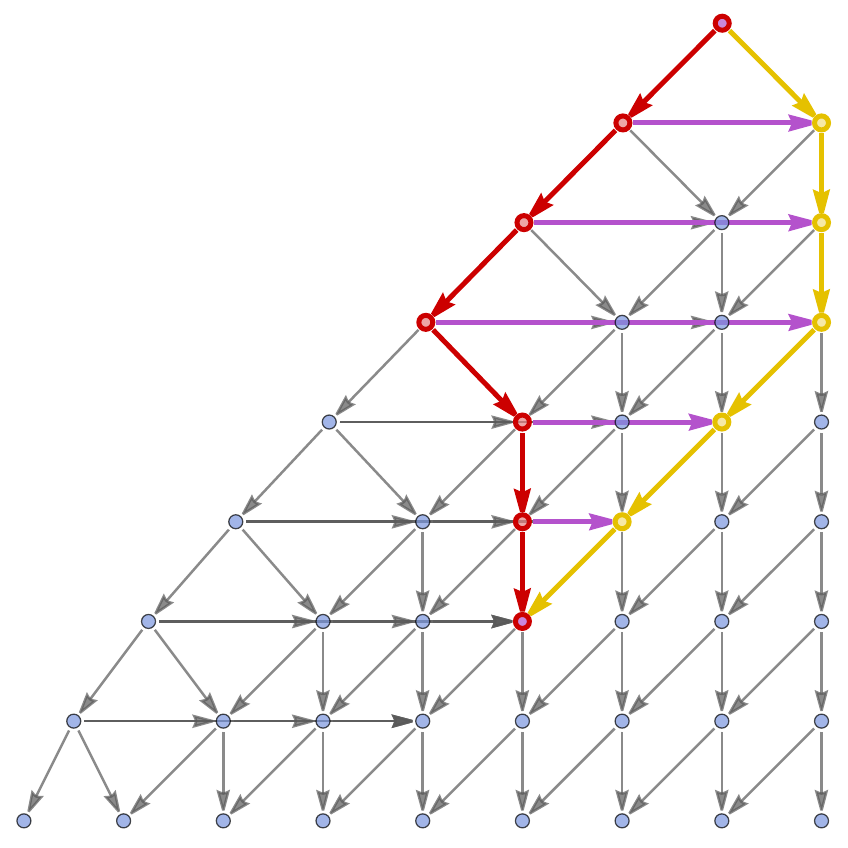}
\caption{On the left we show the construction of Figure \ref{fig:Figure4} (with ad hoc homotopy maps) in flattened-out rendering. On the right is a multiway evolution graph with explicit Wolfram model rules that map vertices on the red path to corresponding vertices on the yellow path. The purple paths still indicate homotopy maps, this time resulting from explicit Wolfram model rules. The two multiway graphs thus constructed are certainly not identical as the new rewriting rules inducing homotopies between the selected paths, are also applicable to other states of the multiway system.}
\label{fig:Figure5}
\end{figure}

\begin{figure}[ht]
\centering
\includegraphics[width=0.5\textwidth]{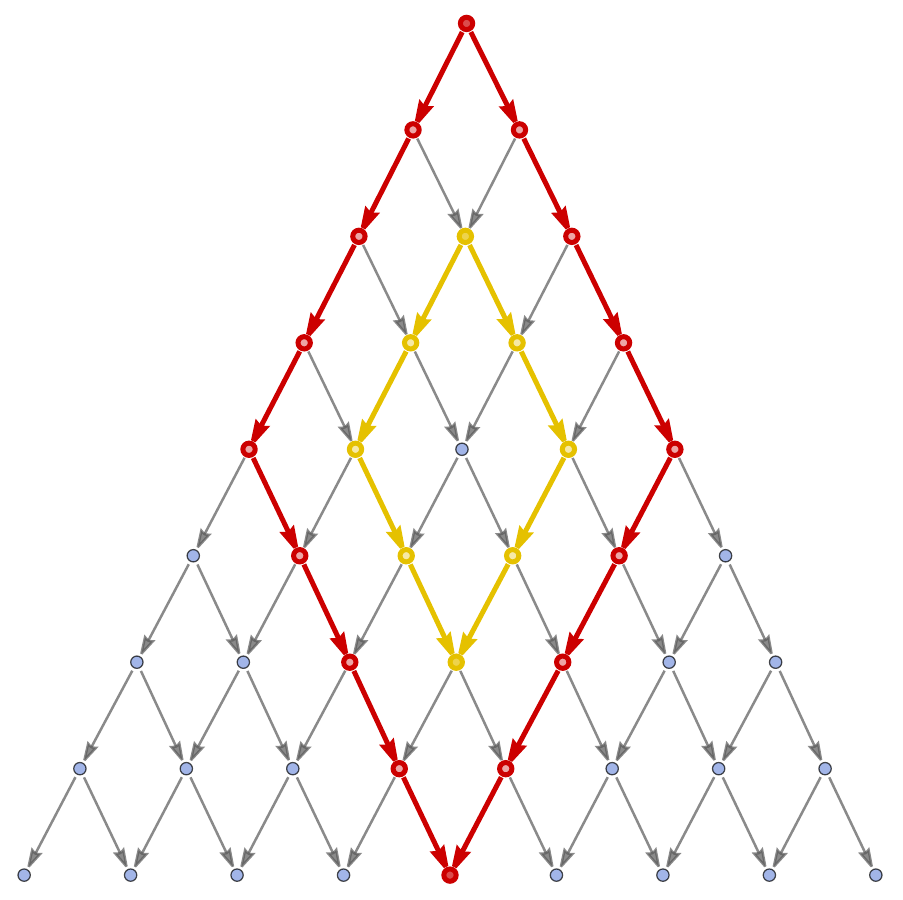}
\caption{In order to illustrate homotopy 3-cells in what follows, we now select four proof paths in our multiway evolution graph. The two red paths indicate two proof paths between vertices ${AA}$ and ${ABBBBABBBB}$, whereas the two yellow paths illustrate two proofs of  ${ABAB \to ABBBABBB}$. }
\label{fig:Figure7.2}
\end{figure}

\begin{figure}[ht]
\centering
\includegraphics[width=0.5\textwidth]{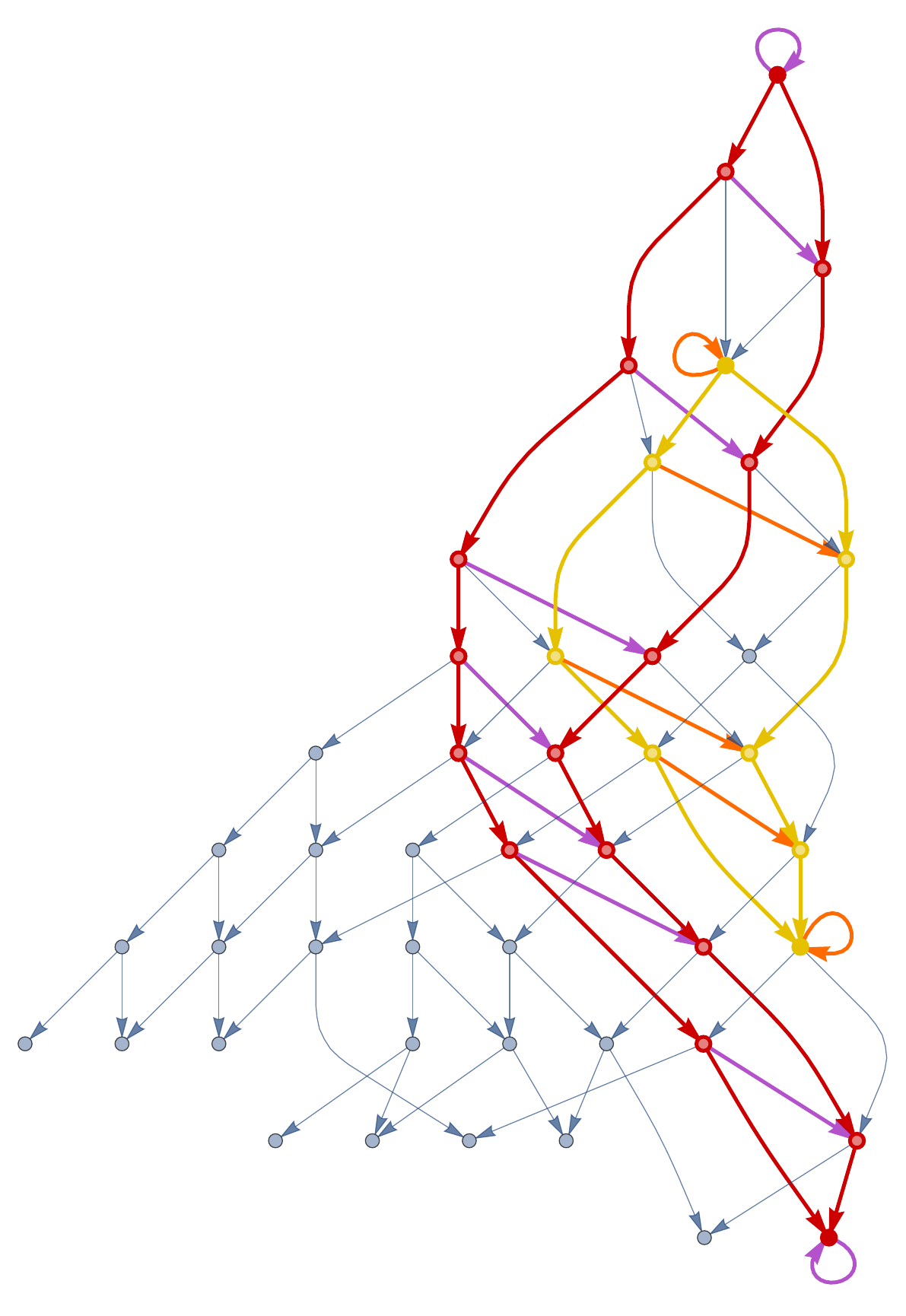}
\caption{A multiway system with homotopy 3-cells. The purple paths create 2-cells between the red paths, whereas the orange paths yield 2-cells between the highlighted yellow paths. The 2-cells here are composed of squares as morphisms of a double category. The lighter arrows between two squares form cubes, yielding morphisms of a 3-fold category. Composition of cubes gives the 3-cell in the figure between the red and yellow 2-cells. Here, all homotopy maps have been externally imposed upon this multiway system. Later, we shall recover these homotopies completely from multiway rewriting rules. }
\label{fig:Figure6}
\end{figure}

This general algorithmic procedure for introducing homotopy maps between paths in an arbitrary multiway system may now be iterated, so as to introduce new homotopy maps between these homotopies. For instance, in Figure \ref{fig:Figure7.2} we make a selection of four path between which we construct 2- and 3-cells as shown in Figure \ref{fig:Figure6} (using the ad hoc method) and then in Figure \ref{fig:Figure7}  (using Wolfram model rewriting rules). 

\begin{figure}[ht]
\centering
\includegraphics[width=1.0\textwidth]{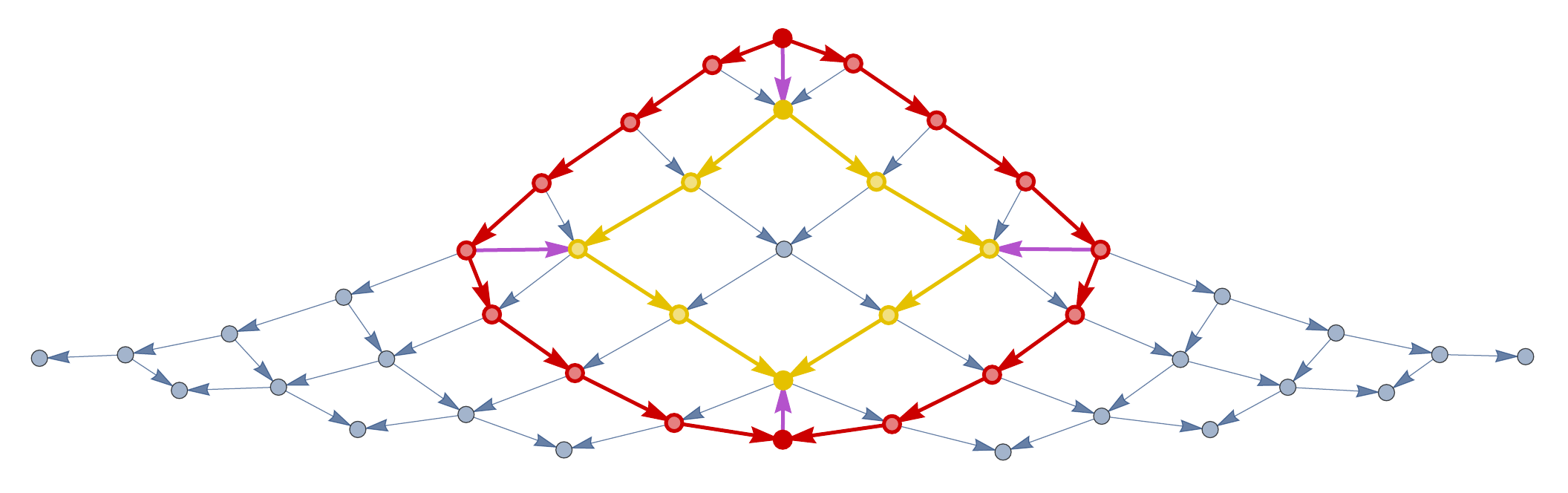}
\caption{An alternate representation of a 3-cell formed by two 2-cells, bounded by red and yellow morphisms respectively, and purple arrows as morphisms between 2-cells in a 3-fold category. }
\label{fig:Figure7.1}
\end{figure}

\begin{figure}[ht]
\centering
\includegraphics[width=0.4\textwidth]{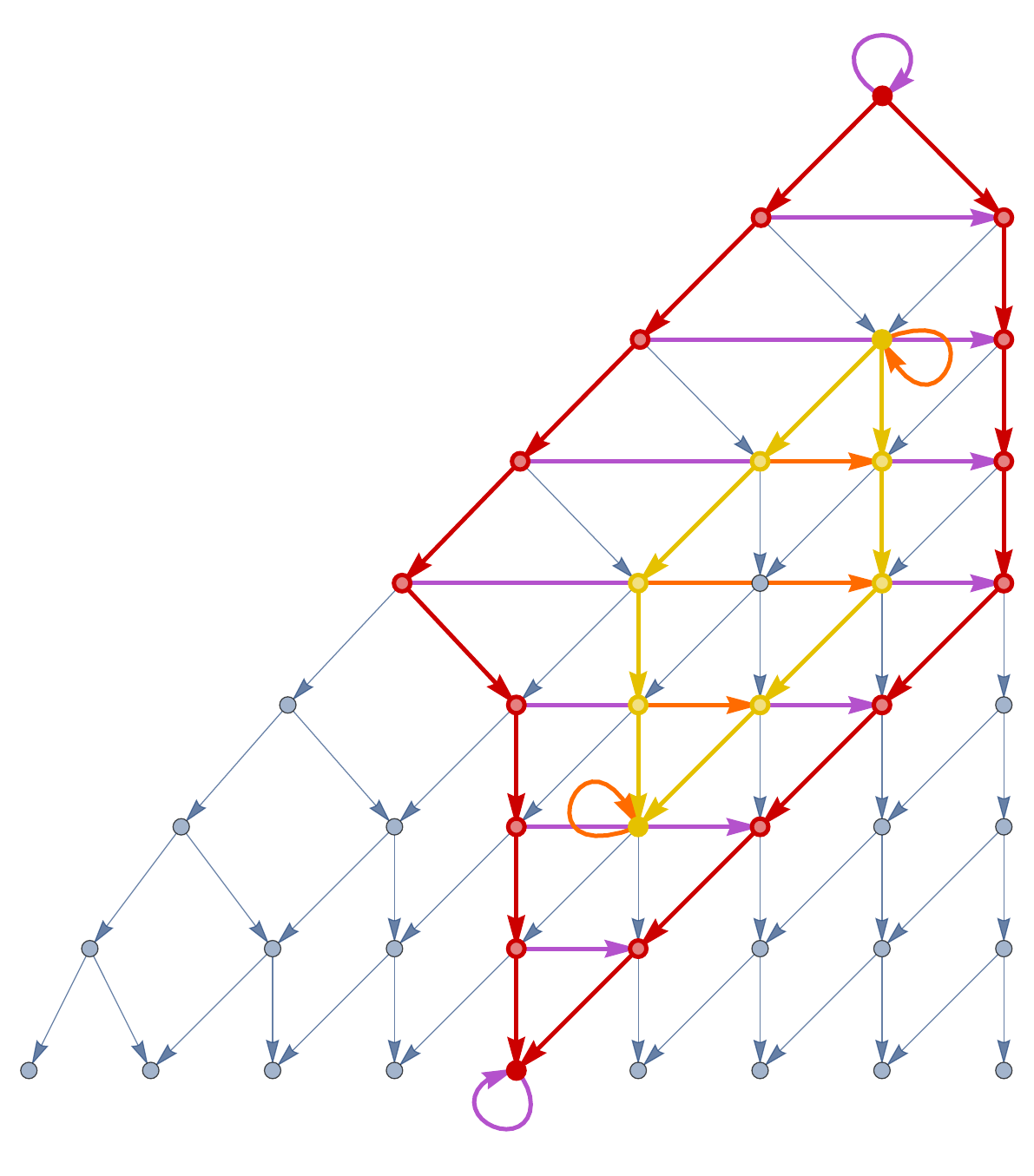}\hspace{0.1\textwidth}
\includegraphics[width=0.4\textwidth]{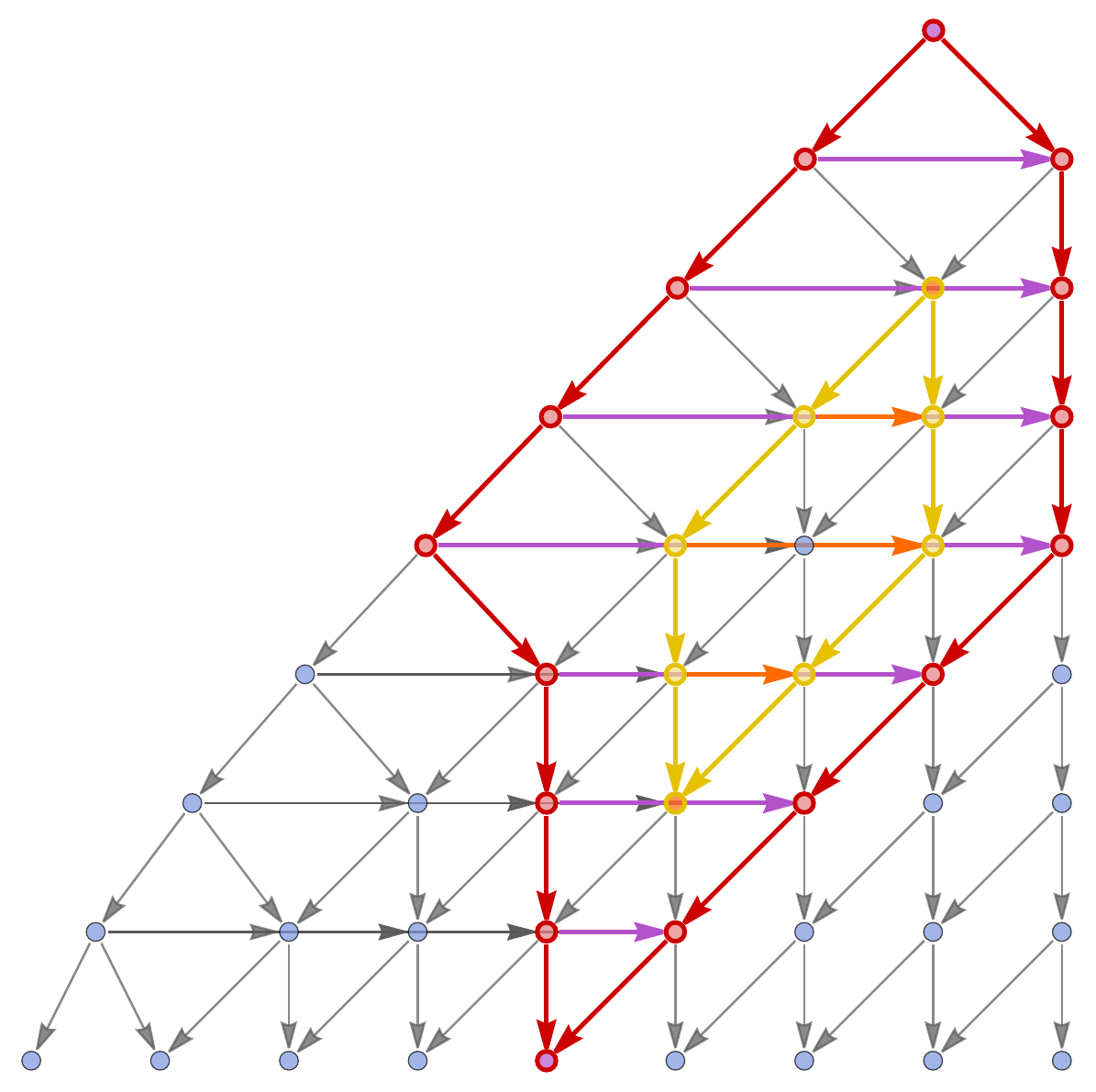}
\caption{The left figure shows the same system as in Figure \ref{fig:Figure6} (with externally imposed homotopy maps) in flattened-out rendering. On the right side, we have a multiway system with explicit Wolfram model rules for constructing the homotopy 2- and 3-cells alluded to in Figure \ref{fig:Figure6}. In addition to inducing higher homotopies, the new rewriting rules also result in new 1-morphisms shown in the multiway evolution graph on the right. }
\label{fig:Figure7}
\end{figure}

Of course, there are many ways of introducing such homotopy maps. An alternate representation of a 3-cell formed by two 2-cells, bounded by 1-morphisms is shown in Figure \ref{fig:Figure7.1}. This construction uses the rules
\begin{multline}
\hspace{2.7cm}  \left\lbrace AA \to ABAB, \,  ABBBBABBBB \to ABBBABBB,  \right.\\  
\left. AABBBB \to ABABBB, \,  ABBBBA \to ABBBAB  \right\rbrace  \hspace{3cm}
\end{multline}
to build arrows between 2-cells.

As mentioned earlier, all these constructions work for any class of multiway system given that the above homotopy maps only really depend on structural properties of these systems and that the collection of homotopy inducing rewriting rules are merely mappings between structures (or cells). 

In the context of homotopy type theory, types equipped with order-$n$ homotopies are referred to as homotopy n-types. Hence, multiway systems equipped with such structures realize homotopy n-types (this is formalized in terms of higher categories in the following subsection). 

\clearpage

\subsection{Higher Categorical Formulation }

We now formalize the above constructions in higher category theory. Starting with the order-2 homotopies between paths (1-morphisms) we make the following proposition: 
\begin{proposition}
The multiway rewriting system shown in Figure \ref{fig:Figure5} is an order-2 homotopy rewriting system (or a homotopy 2-type)  with 2-morphisms between paths,  yielding a double category.
\label{p3.1}
\end{proposition}
Let us first recall the definition of a double category.

\begin{definition}
A double category $D$, denoted  \begin{tikzcd}  D_1 \arrow[r, shift left]   \arrow[r, shift right]  & D_0   \end{tikzcd}  is defined in terms of the following conditions:
\begin{enumerate}[(i)]
\item
The objects of $D$ are the objects of $D_0$ 
\item
$D$ has vertical morphisms, which are the morphisms of ${D_0}$ 
\item
Additionally, $D$ also has horizontal morphisms, which are the objects of ${D_1}$  
\item
Finally, the 2-morphisms of $D$ (also referred to as squares or 2-cells) are the morphisms of ${D_1}$ 
\end{enumerate}
\end{definition}
Here, a 2-cell in $D$ can be represented by the following commutative square:
\begin{equation}
\begin{tikzcd}[column sep={3cm,between origins},row sep={3cm,between origins}]
a \dar{l} \rar["f",""{name=U, below}""] & c \dar{m} \\
b \rar["g" below,""{name=D}""] & d \\
\arrow[Rightarrow, from=U, to=D,"\phi"]
\end{tikzcd} 
\end{equation}
where $a$, $b$, $c$, $d$ are objects; $l$, $m$ are vertical arrows; $f$, $g$ are horizontal arrows; and $\phi$ denotes the 2-cell. 

Note that a double category is an internal category in {\bf Cat}. Vertical composition in a double category is given by composition in the categories $D_0$ and $D_1$, while horizontal composition is given by composition on \begin{tikzcd}  D_1 \arrow[r, shift left]   \arrow[r, shift right]  & D_0   \end{tikzcd} itself, by virtue of it being a category internal to {\bf Cat}. This set-up will be useful for proving the desired proposition.

{\it Proof of Proposition \ref{p3.1}}: 
Following the functorial description provided in section 2, any abstract  multiway system may be thought of as a category ${\cal M}_0$ whose  objects are given by states (rewriting terms) of the multiway system, and whose paths of causally-ordered rewriting chains refer to morphisms of this category. Composition is given by concatenation of directed paths. Transitive closure along directed paths ensures associativity and self-loops on nodes (trivially obtained via identity rewriting rules) ensure the identity axioms of a category.

From the definition of a double category given above, ${\cal M}_0$ is precisely the category $D_0$, with objects being the vertices of a multiway evolution graph (a prototypical example being given in Fig. \ref{f3.1.1}),  and paths formed by multiway evolution (also shown in Figure \ref{fig:Figure3}) corresponding to vertical morphisms in $D_0$. 
The horizontal arrows in Figure \ref{fig:Figure5}, colored in purple, are induced by the introduction of additional rules to the rewriting system and serve the purpose of connecting (vertical) paths in the multiway system. Hence, these can be identified as objects of another category, which we denote as ${\cal M}_1$  (${\cal M}_1$ will be identified with $D_1$ once we specify its morphisms).  
Imposing commutativity over squares formed by the vertical and horizontal arrows of the multiway system gives us the morphisms of ${\cal M}_1$; these squares, exemplified in Figure \ref{fig:Figure5}, are 2-morphisms of the double category   \begin{tikzcd}  {\cal M}_1 \arrow[r, shift left]   \arrow[r, shift right]  & {\cal M}_0   \end{tikzcd}. 

Given this double category, the 2-morphism in the multiway evolution graph between the highlighted red and yellow paths in Figure \ref{fig:Figure5} is precisely given by vertical compositions of squares. These are the 2-cells of the new multiway rewriting system, enhanced with additional rules, thus making it an  \textit{order-2 homotopy rewriting system} or a homotopy 2-type. (Note that the triangles at the start and end points of the 1-paths are also squares with an identity morphism at those points.)
\qed \\  \; 

This leads us directly to the following corollary: 
\begin{corollary}
\label{c3.1}
A  homotopy 2-type multiway system extends to a double groupoid upon admitting invertible rules to those used in the multiway construction (including the homotopy rules).
\end{corollary}

{\it Proof of Corollary \ref{c3.1}}: 
Adding invertible rules to the multiway system equips every causal arrow with an inverse arrow, such that their composition is the identity map on either rewriting term. Given a homotopy 2-type multiway system, the presence of invertible arrows implies that all 1-morphisms and 2-morphisms are isomorphisms. This extends the multiway system in Proposition \ref{p3.1} to a double groupoid.  \qed 

Furthermore, as shown in the previous section, with appropriate additional rules, one can continue the above construction inductively to obtain rewriting systems with progressively higher homotopies between paths. This yields the following proposition:

\begin{proposition}
So long as additional rewriting rules between cells, up to order $n - 1$, are admissible, then the multiway rewriting system shown in Figure \ref{fig:Figure5} can be enhanced to yield a homotopy n-type  with $n$-morphisms between paths, thus yielding an $n$-fold category.
\label{p3.2}
\end{proposition}

{\it Proof of Proposition \ref{p3.2}}: 
The proof of this proposition follows by a simple inductive construction. 

A double category itself is a  special case of an $n$-fold category. In Proposition \ref{p3.2}, we proved this proposition for the $n = 2$ case.  

For the $n = 3$ case, we need to construct an order-3 homotopy rewriting system equipped with 3-morphisms, which are cubes bounded by squares of the given double category. This can be done by adding specific rewriting rules that give arrows between paths (1-cells) as shown in    Figures \ref{fig:Figure6} and \ref{fig:Figure7}  such that one obtains multiple distinct 2-cells, which eventually permit additional arrows between them, giving 3-morphisms as cubes bounded by squares. These cubes define a 3-fold category as follows:
\begin{definition}
A 3-fold category, denoted   \begin{tikzcd}  D_2 \arrow[r, shift left]   \arrow[r, shift right]  & D_1 \arrow[r, shift left]   \arrow[r, shift right]  & D_0   \end{tikzcd}   is defined in terms of the following conditions:

\begin{enumerate}[(i)]
\item
The objects are the objects of $D_0$ 
\item 
The vertical arrows are the morphisms of $D_0$ 
\item
The horizontal arrows are the objects of $D_1$  
\item
The vertical squares are the morphisms of $D_1$ 
\item
The horizontal squares are the objects of $D_2$  
\item
The cube bounded by vertical and horizontal squares is a morphism of $D_2$ 
\end{enumerate}
The above data has to satisfy the commutative diagram:
\begin{equation}
\begin{tikzcd}[row sep=1.5em]
A \arrow[rr,"f",""{name=0, below}] \arrow[dr,swap,"a",""{name=1}""] \arrow[dd,swap,"h",""{name=2}""] &&
  B \arrow[dd,swap,"h'" near end,""{name=3}""] \arrow[dr,"b",""{name=4, below}""] \\
& A'\arrow[rr, crossing over, "f'" left,""{name=U, below}"", yshift=-0.05ex]
&& B' \arrow[dd,"k'",""{name=5}""] \\
C \arrow[rr,"g" near start,""{name=N, below}""] \arrow[dr,swap,"c",""{name=6}"", below] && D \arrow[dr,"d",""{name=7, below}""] \\
& C' \arrow[rr,"g'" near start,""{name=M}""] \arrow[uu,<-,crossing over,"k" near end]&& D'
 \arrow[from=1,to=4, dash, shorten=10mm, phantom,""{name=K}""]
 \arrow[r,from=6,to=7, dash, shorten=10mm, phantom,""{name=L}""]
 \arrow[Rightarrow, line width=0.3pt, from=K, to=L,"\psi"] \\
 \arrow[r, dash, line width=0.3pt, from=K, to=L]
 \arrow[Rightarrow, from=0, to=N, xshift=2ex, crossing over, shorten=1mm, "\eta", near start, swap]
 \arrow[Rightarrow, from=U, to=M, xshift=2ex, crossing over, shorten=1mm, "\omega", near end, swap]
 \arrow[Rightarrow, from=1, to=6, shorten=2mm, "\alpha",near start]
 \arrow[Rightarrow, from=4, to=7, shorten=1mm, "\beta"]
 \arrow[r, from=1, to=4, shorten =11mm, xshift=10ex, yshift=-3.03ex, crossing over]
 \arrow[r, dash, from=1, to=4, shorten=11mm, xshift=-0.3ex, yshift=-3.03ex, crossing over]
 \end{tikzcd} 
\end{equation}
where arrows are 1-morphisms, squares are 2-morphisms and cubes are 3-morphisms (3-cells) of the 3-fold category.   Compositions, identity and associative laws on each of these entities follow from their usual definitions in $D_0$, $D_1$ and $D_2$ where appropriately applicable.       
\end{definition}
The 3-cells in our multiway rewriting system in Figure \ref{fig:Figure7}  are then given by vertical, horizontal and sideways compositions of cubes of a 3-fold category.

Likewise, one can continue this process of constructing higher morphisms indefinitely, so long as additional rewriting rules between cells of up to order $n - 1$ are admissible such that $n$-morphisms are now $n$-hypercubes defined within an $n$-fold category, and $n$-cells of such multiway systems are simply compositions of $n$-hypercubes in $n$-directions. The ensuing  $n$-fold category can then be expressed as:  
\[  \begin{tikzcd} {\cal M}_{n - 1} \arrow[r, shift left] \arrow[r, shift right]  & {\cal M}_{n - 2} \,\,  \cdots\cdots \,\, {\cal M}_2 \arrow[r, shift left] \arrow[r, shift right] & {\cal M}_1 \arrow[r, shift left] \arrow[r, shift right]  & {\cal M}_0 \end{tikzcd} \]  
where the objects of ${\cal M}_i$ (for $0 < i < {n-2}$) are $i$-dimensional hypercubes in ${\mathbb R}^n$ whose normal vectors are oriented along the $z$-axis in ${\mathbb R}^n$. The morphisms of ${\cal M}_i$ are $(i + 1)$-dimensional hypercubes in ${\mathbb R}^n$ whose normal vectors are oriented orthogonal to the $z$-axis in ${\mathbb R}^n$. Furthermore, the objects of ${\cal M}_{n - 1}$  are $(n - 1)$-dimensional hypercubes with normal vectors oriented along the $z$-axis in ${\mathbb R}^n$ and the morphism in ${\cal M}_{n - 1}$ is the commutative $n$-hypercube composed of the above. 

This iterative definition of an $n$-fold category corroborates the previously described construction, in which homotopies of up to order $n$ in multiway systems are realized by including additional rewriting rules that introduce new arrows to obtain squares (and so on), until one eventually obtains a commutative $n$-hypercube. 
\qed  \\  \,

The benefit gained by realizing $n$-fold categories in our constructions is the ease of expressing multiple compositions by gluing hypercubes up to order $n$. Composable arrays of hypercubes can then be used to construct all ($i \leq n$)-cells of the category. 

For example, in a 2-fold category, we can define a composable array of 2-dimensional elements (squares) to be such that any array is composable with its immediate neighbors (the associative and commutative laws imply that the composition is well defined). This process easily extends in $n$-fold categories  using $n$-dimensional elements ($n$-hypercubes). 

A few remarks are in order:
\begin{remark}
As in Corollary \ref{c3.1}, a homotopy n-type multiway system extends to an $n$-fold groupoid upon admitting invertible rules that ensure invertibility of all higher morphisms.
\label{rem:Remark1}
\end{remark}

\begin{remark}
Even though an $n$-fold category is a strict version of an $n$-category, in that all $n$ composition operations are strictly unital and associative and strictly commute with each other, nonetheless, $n$-fold groupoids model all    homotopy n-types, and this fact is what ensures that multiway rewriting systems are models of  homotopy n-types.
\end{remark}

\begin{remark}
Note that the propositions above have been described referring to string substitution multiway systems only for ease of illustration. However, it is easy to see that these statements hold for any class of multiway system as the proofs above are only based on structural properties of these systems and do not depend on the specific type of rewriting states used.
\end{remark}

\subsection{The $\infty$-Limit of Multiway Rewriting Systems as $\infty$-Groupoids}

Given the iterative construction presented above of inducing higher homotopies in multiway rewriting systems via the inclusion of supplementary rules, one can now ask what the $n \to \infty$ limit of this construction yields? 

\begin{proposition}
The $n \to \infty$ limit of a  multiway rewriting system with invertible homotopy rules, corresponding to higher morphisms (so long as they are admissible), is an $\infty$-groupoid. 
\label{p3.3}
\end{proposition}

{\it Proof of Proposition \ref{p3.3}}: 
Following Proposition \ref{p3.2} and Remark \ref{rem:Remark1}, multiway rewriting systems equipped with homotopies up to order $n$, along with invertible rules, yield an $n$-fold groupoid. 
The $n \to \infty$ limit of an $n$-fold groupoid is precisely an $\infty$-groupoid. This limit exists so long as the infinite hierarchy of rules required to iteratively construct higher morphisms is admissible on that multiway rewriting system.  
\qed \\ \, 

The following remark is due:
\begin{remark}
The $\infty$-groupoids discussed here have been constructed within the model of strict n-categories. Even though some definitions of $\infty$-categories (or statements about the homotopy hypothesis based on those definitions) that declare  $\infty$-groupoids as topological spaces are formulated within the model of weak or quasi categories, there is indication that such definitions or statements will eventually be expressible within other models of categories   \cite{riehl20152},  \cite{Riehl2017a},  \cite{riehl2018elements}.  Hence, it would be interesting to consider what such categorical weakenings would be interpretable as within the kinds of rewriting systems considered above. 
\end{remark}

To sum up, upon expressing Wolfram model multiway systems as homotopy types and associated multiway rules as type constructors, what we find here is that inducing a homotopy between two paths (in other words, applying a completion procedure between multiway branches) can be thought of as introducing new higher-order  rewriting rules, which have the effect of producing higher-order structures from the one we originally started with. This notion of higher-order structures is  made more precise using the language of category theory; if each path is interpreted as a morphism between objects (known as a 1-morphism), then a homotopy can be interpreted as a morphism between 1-morphisms (that is, a 2-morphism), with the resultant structure being a 2-fold category (depending on the definition of higher categories one uses).  Homotopies of 2-morphisms can then be interpreted as 3-morphisms within 3-fold categories, and so on, thus producing a whole infinite hierarchy of higher-order categories. The limit of this hierarchy is the $\infty$-groupoid, which can be thought of as being the structure obtained by inducing all possible homotopies (and, consequently, by applying all possible completion procedures) to a given multiway system.  On the other hand, replacing isomorphisms with strict equalities would once again collapse this tower of homotopies to obtain a 1-category. In that sense, an $\infty$-groupoid is a natural structure with fewer constraints required to truncate higher morphisms.

Interestingly, we can generalize the above structures even further in higher categories. For this, we will need the following definitions:

\begin{definition}
\textbf{Rulial Multiway System}: We define the `Rulial Multiway System' as the multiway system equipped with higher homotopies up to order $n > 1$. In other words, these systems are homotopy n-types. Furthermore, the limit of this structure when $n \to \infty$ will be referred to as the `Limiting Rulial Multiway System'.    
\end{definition}

\begin{definition}
\textbf{Rulial Multiverse}: We define the `Rulial Multiverse' as the collection of all possible multiway systems of a given class, including all rulial multiway systems (of that class). Further, the specific collection of all limiting rulial multiway systems of a given class will be called the limiting rulial multiverse.  
\label{def3.4}
\end{definition}

This leads to the following proposition: 
\begin{proposition}
The limiting rulial multiverse, wherein all limiting rulial multiway systems are equipped with invertible homotopy rules,  yields  an $\left(\infty, 1\right)$-category  {\bf $\infty$Grpd} of $\infty$-groupoids. 
\label{p3.4}
\end{proposition}

{\it Proof of proposition \ref{p3.4}}: 
An $\infty$-groupoid is an $\left(\infty, 0\right)$-category, where all morphisms (up to order $\infty$) are invertible. We denote this as $Cat_{\left(\infty, 0\right)}$. In proposition \ref{p3.3} above, we have shown that the $\left(n \to \infty \right)$-rulial multiway rewriting system with invertible homotopy rules is a type-theoretic construction that realizes $Cat_{\left(\infty, 0\right)}$. By definition \ref{def3.4}, these $\infty$-categories are among the objects of the full rulial multiverse, which thus internally realizes {\bf $\infty$Grpd}, the  $\infty$-category of $\infty$-groupoids. Furthermore, from definition 1.3.6, remark 1.3.7 and corollary 4.3.16 in Lurie \cite{lurie2009infinity}  it follows that {\bf $\infty$Grpd} is equivalent to the $\infty$-category $CSS_{Cat_{\left(\infty, 0\right)}}$ of complete Segal space objects comprising the $\infty$-category $Cat_{\left(\infty, 0\right)}$. And from  \cite{lurie2009infinity},  \cite{lurie2009higher}  this construction of a complete Segal space models an $\left(\infty, 1\right)$-category. This $\left(\infty, 1\right)$-category is internal to the full rulial multiverse. Therefore, the limiting rulial multiverse yields an $\left(\infty, 1\right)$-category.  
\qed

\begin{remark}
Using definition 1.3.6 of Lurie  \cite{lurie2009infinity}  an $\left(\infty, n\right)$-category $Cat_{\left(\infty, n\right)}$ can be defined inductively as the $\infty$-category $CSS_{Cat_{\left(\infty, n-1\right)}}$ of complete Segal space objects comprising the $\infty$-category $Cat_{\left(\infty, n-1\right)}$. This suggests a hierarchy of inductive generalizations of the limiting rulial multiverse itself. 
\end{remark}

Starting with an $Cat_{\left(\infty, 2\right)}$ category that is internal to the multiverse of multiverses, which we denote as a 2-fold multiverse, one can formally express n-fold multiverses. However, beyond the 2-fold case, it is not clear what the precise physical interpretation of these structures would be. The 2-fold multiverse corresponds to a statistical ensemble of limiting rulial multiverses where each multiverse represents a distinct fibration over rulial space. Hence, the choice of a global geometry for the limiting rulial multiverse can be thought of as the analogue of  fixing a gauge in the ensemble of multiverses.

\begin{remark}
The {\bf $\infty$Grpd} is often used as the archetypical example to define an  $\left(\infty, 1\right)$-topos, the home of classical homotopy theory.  In  \cite{Schreiber2013a},  it was shown that an $\left(\infty, 1\right)$-topos of sheaves ${\bf Sh_{\infty} (C)}$ over a classifying  $\infty$-category $C$, when equipped with the Grothendieck topology, comes equipped with a triplet of adjoint functors to  {\bf $\infty$Grpd}, which preserve discrete and indiscrete topology. Furthermore, if the fundamental $\infty$-groupoid functor also exists and is adjunct to this triplet, that preserves a cohesive structure, which synthetically defines formal geometric spaces (as smooth sets or continuous sets) within ${\bf Sh_{\infty} (C)}$.  This corroborates with our constructions in this work   wherein the limiting rulial multiverse provides the underlying computational model (as a rewriting system) of formal spaces (as $\infty$-groupoids) on which physics can be done.
\end{remark}

Following the above remark, note that the existence (or lack of it thereof) and  construction of cohesive structures on limiting rulial multiway systems will certainly depend on the rules used to construct that multiway graph, which in turn determine the growth of that multiway system (see  \cite{zeschke2021growth}  for a classification on multiway growth functions).    In  \cite{Schreiber2013a}, \cite{Shulman2017}  two concrete examples of cohesive structures have  been reported: those constructed within a topos of convergent sequences or what  are referred to as  `consequential spaces'  $\{ \mathbb{N}_{\infty} \}$; and those constructed within a category of abstract co-ordinate charts $\{ \mathbb{R}^n \}$ (as types, not as sets or topological spaces) for all $n \in {\mathbb N}$.  In future work, it will be interesting to investigate other realizations of cohesive structures, using multiway systems.

Another issue that we have not discussed very much about here concerns local presentations of the above $\infty$-categories. Some of these can be expressed in terms of simplicial sets and Kan complexes. For instance, in classical homotopy theory, an {\bf $\infty$Grpd} is the simplicial localization of the category {\bf sSet} of simplicial sets. It can also be presented as a Kan-complex enriched category. Namely, it is the full {\bf sSet}  enriched subcategory of {\bf sSet} on the Kan complexes. These local presentations of  $\infty$-categories and their multiway rewriting models may presumably be   relevant to the study of local operator algebras, potentially in the context of axiomatic quantum field theories, such as those discussed in    \cite{schreiber2009aqft}.

\section{Pregeometric Structures Internalized within an $\infty$-Topos}

As we have just seen, the Wolfram model multiway rewriting systems above are higher categorical objects constructed as a type theory. These were all based on purely syntactic or combinatorial definitions that did not hinge upon any a priori notion of space or geometry. In other words, these are true pregeometric structures (as opposed to those that are defined as discretization of an underlying geometry). Hence, an obvious challenge for this model is to answer how spatiality or geometry may even emerge from these constructions.   This is what we address in this section. The framework of homotopy type theory ties closely with current developments in higher category theory. Following Grothendieck's Homotopy Hypothesis  (see \cite{baez2007homotopy} for a review),  $\infty$-groupoids are identified as formal spaces. Moreover, our pregeometric constructions above can be internalized within a suitable $\infty$-topos, within which, spaces and constructions relevant to physics may be realized (as is also being attempted in  \cite{Schreiber2013a}).

\subsection{A Perspective on the Homotopy Hypothesis and Synthetic Geometry}


As we've seen, a multiway rewriting system can be endowed with higher homotopies  by  systematic application of additional rewriting rules at each homotopy level, resulting in paths between paths and so on.  The limiting $n \to \infty$ structure one thus obtains is the limiting rulial multiway system. Additionally, admitting inverse morphisms  (which can also be thought of as enacting a certain kind of localization), this limiting  structure can be identified as  an  $\infty$-groupoid. By itself, this realization is significant because of Grothendieck's Homotopy Hypothesis, which relates $\infty$-groupoids to  formal (homotopy) spaces (which, within certain models of higher category theory, are identified as topological spaces \cite{baez2007homotopy}).  If Grothendieck's hypothesis is true, it then provides rigorous justification for the claim that the limiting rulial multiway system is itself a topological space, and    may be endowed with a non-trivial topology. In our model, this provides an entry point for synthetic topology (and subsequently synthetic geometry), where we started with purely combinatorial structures as building blocks of a type theory and now functorially realizes potential topological (and with additional conditions, possibly also geometrical) spaces as limiting structures within a classifying category.  Moreover, as has been argued in  \cite{Shulman2016},  \cite{Shulman2017}   within the broader context of homotopy type theory, one can make this correspondence from type theory to spaces even more specific. We now describe the implications of that specification for limiting rulial multiway systems in our model. 

At its core, synthetic geometry via homotopy type theory, championed by     \cite{Shulman2017}   seeks to investigate how non-trivial  spatial structures (topological or geometric) are associated to a variety of algebraic objects such as groups, rings, lattices, etc., almost all of which can be realized purely from type-theoretic foundations. To do so, one then needs to determine how these syntactic definitions can be internalized within a suitable topos (or higher topos). Furthermore, in order that a classifying category of a type theory be elevated to the status of a topos, one requires additional constructions within the syntax category of the type theory. Firstly, one needs to include as type constructors a subobject classifier $\Omega$ (which is the categorical generalization of a subset identifier) as well as a `universe type' ${\cal U}$.  One also requires finite limits and colimits to exist (including generalizing to homotopical products and pullbacks when working with higher categories). In fact, once one equips a Wolfram model multiway system with these additional constructions, the resulting system can be interpreted within a suitable elementary free topos  ${\cal T}$ (also referred to as the logical topos).  Following  \cite{Shulman2017}, one can then  map  ${\cal T}$ to a  topos of spaces  ${\cal S}p$ and ask whether type constructions in ${\cal T}$ faithfully preserve spatial structures contained in ${\cal S}p$?  To address this question, one investigates the following adjunctions: 

\begin{eqnarray}
\begin{tikzcd}[column sep=7em, row sep=6em]  
{\cal T} \,  \arrow[dr, shorten=0.2em, "unique" ',  shift right=3.5]     \arrow[r, "unique",  shift left=0]  & \,\,\,\,  {\cal S}p \,\,  \arrow[d,  "\Gamma" ]   \arrow[d, hookleftarrow,  "\nabla", shift left=6]   \arrow[d, hookleftarrow,  "\Delta",  shift right=6]   \\ 
  &\,\,  {\bf Set}  
\end{tikzcd} 
\label{eq5.1.1}
\end{eqnarray}

Here ${\cal S}p$  comes naturally equipped with a forgetful functor $\Gamma$ to the topos of sets ${\bf Set}$ along with left and right adjoint functors $\Delta$ and $\nabla$ which respectively preserve the discrete and indiscrete topology in ${\cal S}p$. As proposed in  \cite{Shulman2017},  so long as the free topos (associated to the underlying type theory) with these string of adjunctions \textit{factors uniquely} through the topos of spaces, the corresponding  type theory is endowed with spatial structure from  ${\cal S}p$.  Note however, that the adjunctions in eq. (\ref{eq5.1.1})  do not preserve infinite limits or function-spaces with non-discrete domains. In other words, it only guarantees discrete topological spaces for the associated type constructions. 

The question then is: How does one extend the above formalism to enable type constructions to synthetically inherit non-discrete topologies and smooth spaces? This was addressed in \cite{Schreiber2012}, \cite{Schreiber2013a}    using what is called a cohesive $\infty$-topos (see also  \cite{lawvere2003sets}   for the axiomatization of cohesive toposes).  A  cohesive $\infty$-topos is one that is local and  $\infty$-connected. In this context, a `cohesive structure' replaces the category ${\bf Set}$ above with the category  ${\bf \infty Grpd}$. In that, higher homotopies and $\infty$-groupoids naturally  enter this picture. The functorial construction involves an additional left adjoint   to the triplet of adjunctions in eq. (\ref{eq5.1.1}), but now between the topos of 
sheaves ${\bf Sh_{\infty} (C)}$ and ${\bf \infty Grpd}$ (shown in eq.   (\ref{eq5.1.2})).  This left adjoint is the  fundamental $\infty$-groupoid functor $\Pi_{\infty}$.  ${\bf Sh_{\infty} (C)}$ is an  $\left(\infty, 1\right)$-topos of sheaves  over a classifying  $\infty$-category $C$, and comes equipped with the Grothendieck topology.    The existence of  $\Pi_{\infty}$  (also referred to as the ``shape functor")  ensures  cohesive structure (locality and  $\infty$-connectedness of  ${\bf Sh_{\infty} (C)}$),  and the corresponding homotopy type theories in $C$ then synthetically inherit topological or smooth structures (referred to as `continuous sets' or `smooth sets' respectively in    \cite{Schreiber2013a}).  Just as the ordinary fundamental groupoid functor ensures that topological data up to order-1 homotopy is preserved, the  existence of the fundamental $\infty$-groupoid  functor now ensures that topological data at all orders of homotopy are preserved.  Eq. (\ref{eq5.1.2})  below shows the relevant diagram of adjunctions for constructing  cohesive structures: 

\begin{eqnarray}
\begin{tikzcd}[column sep=7em]  
\infty Func\left( C^{op}, {\bf \infty Grpd} \right)  \arrow[r, hookleftarrow, shift right=2]   \arrow[r, "\infty-Stackification",  shift left=2]  & {\bf Sh_{\infty}}(C)  \arrow[r,  "\Gamma", shift right=2.5 ]  \arrow[r,  "\Pi_{\infty}",  shift left=7.5]   \arrow[r, hookleftarrow,  "\Delta", shift left=2.5]   \arrow[r, hookleftarrow,  "\nabla",  shift right=7.5]  & {\bf \infty Grpd}
\end{tikzcd} 
\label{eq5.1.2}
\end{eqnarray}  
  
Certainly, not every category $C$ will admit cohesive structures or realize the construction in eq. (\ref{eq5.1.2}).  Among the ones that do, are the category of abstract co-ordinate charts and the category of consequential spaces  \cite{Schreiber2013a}, \cite{Shulman2017}. Hence, it would be useful to find additional realizations of  cohesive structures from type constructions, particularly those from generic multiway systems.


\subsection{A Comment on the Rulial Multiverse as a Fibration }

In the previous section, we have argued that the limiting Wolfram model rulial multiverse of $\left(n \to \infty \right)$-rulial multiway rewriting systems with invertible homotopy rules limits to an $\left(\infty, 1\right)$ category  {\bf $\infty$Grpd} of $\infty$-groupoids.  This {\bf $\infty$Grpd} can be used to define an  $\left(\infty, 1\right)$-topos. If we interpret this total space as a fibration of $\infty$-categories, individual rulial multiway systems correspond to fibers of this topos  (presumably this construction yields an $\left(\infty, 1\right)$-Grothendieck fibration; this issue will be formally explored in an upcoming work). Each fiber is generated from a distinct collection of Wolfram model rules (non necessarily without intersections), yielding, via Grothendieck's hypothesis, a collection of topological spaces in the $\left(\infty, 1\right)$-topos. Furthermore, following \cite{Schreiber2013a},  additional functorial constraints on this topos involving the fundamental $\infty$-groupoid functor, preserve `cohesive structures'  that admits geometry on the objects of the $\left(\infty, 1\right)$-topos. 

A noteworthy feature of the synthetic geometry approach, that we have used in our constructions, is that now one does not require to make ad hoc assignments or assumptions of geometric structures on local entities of the model (such as assumptions of a background geometric space or discretization schemes defined using pre-assigned geometric data). Instead, geometry (in the form of cohesive structures) is inherited functorially using global properties,  and as a consequence, is naturally induced upon local structures by taking sections or projections of the total space. For example, foliations are section of the multiway system analogous to spacelike hypersurfaces. These are the so-called  `branchial graphs'  of the Wolfram model, which refer to a network of entangled states. The geometry of these branchial spaces is then the sectional geometry induced via the globular geometry of the rulial multiway system. Another example are single-way paths of the  rulial multiway system. These correspond to classical timelike hypersurfaces (or Wolfram model causal graphs). The induced geometry on these spaces is one obtained via restriction maps to individual multiway paths or homotopy 1-types. Thus, Grothendieck's hypothesis effectively explains why various Wolfram model combinatorial structures, such as multiway systems, branchial graphs, causal graphs, etc., are can be endowed with the spatial structures.


\subsection{Role of Voevodsky's Univalence Axiom in the Topos of Rulial Multiway Systems}

In the foundations of mathematics, the notion of extensionality commonly refers to those criteria by which two objects are deemed identical. For example, the axiom of extensionality in axiomatic set theory states that two sets are identical if and only if they contain the same elements. Extensionality is thus the logical analog of  `state equivalence' in the case of the Wolfram multiway system formalism. However, in conventional mathematical logic, there are several notions of extensionality. For example, propositional extensionality asserts that a pair of propositions may be considered identical if and only if they logically imply each other. Analogously, functional extensionality asserts criteria for identity between two functions. Furthermore, `uniqueness of identity proofs' is yet another extensionality axiom. HoTT attempts to bring these different notions of equality within a common framework. It does so via the univalence axiom, which serves as a grand generalization of and therefore subsumes all other extensionality axioms to what may be called typal extensionality   \cite{Program2013},  \cite{Ahrens2021}.  Namely, that   two types themselves, thought of as $\infty$-groupoids, can be considered identical if and only if they are homotopy equivalent (at all homotopy levels). Interestingly, it turns out that many type theories originating from a constructivist paradigm of mathematics are formal systems that are distinct from ZFC set theory, and do not always admit the law of the excluded middle or the axiom of choice. Rather, many of these type theories and their associated  toposes are consistent with the univalence axiom. All of this has important implications for the very foundations of mathematics. In particular, as has been discussed  in  \cite{Shulman2016}, \cite{Shulman2017},  this suggests that there does not exist one a priori preferred axiomatization to describe mathematics, but that one can in fact  formalize many distinct universes of mathematics based on different axiomatizations and then transform from one formal system to another. 

How does all this relate to Wolfram model multiway systems? When one induces a homotopy between two independent paths in a multiway system (which could also be thought of as a completion in the context of automated theorem proving, or for our purposes here, a proof of equivalence between proofs in the context of homotopy type theory), we effectively treat the corresponding paths as being identical, in the sense that they proceed to evolve as an effective single path in the multiway system. Therefore, implementation of homotopies  or  completion procedures in  multiway systems correspond precisely to propositional extensionality in type theory. This can subsequently be generalized to proof completions in the multiverse itself, that is, to identification of two given multiway system types themselves. The latter is elegantly achieved using the univalence axiom, which generalizes propositional extensionality to type extensionality in the topos of rulial multiway systems. The univalence axiom thus provides the formal backing for the operations of equating states and paths within the combinatorial constructions of the Wolfram model.

\section{Applications to Physics }

The underlying philosophy of the Wolfram model originates from the idea that  the building blocks of the universe are fundamentally discrete entities (and their relations), governed by computational rules. Structures relevant to physics are  constructed from the combinatorics of non-deterministic rewriting systems. A key objective of this current work was to provide a formal foundation for the Wolfram model in homotopy type theory. Doing so, showed how geometry and space emerges from pregeometric structures. We reckon that this framework of higher categorical pregeometric structures may also be relevant to several other approaches addressing questions at the foundations of physics  \cite{isham2000some}, \cite{Schreiber2013a}, \cite{schreiber2016higher}.   Below we identify three specific applications of our formalism.

\subsection{Synthetic Spaces and Geometry from Rewriting Systems using   Homotopy Type Theory}

One of the key takeaways from the synthetic geometry and homotopy type theory program is that the notion of space arises functorially, when the topos of spaces carries adjunctions to the topos of sets, including the existence of the  fundamental groupoid functor.  This means no background geometric spaces or their discretization per se need be assumed. Instead, geometry (in the form of cohesive structures) is inherited by higher structures, and consequently, is induced upon local structures by taking sections or projections of the total space. In the context of homotopy type theory, this topos of spaces is elevated to an  $\left(\infty, 1\right)$-topos, whose objects,  the $\infty$-groupoids, are identified as formal topological spaces (via Grothendieck's hypothesis), and geometric spaces (via additional cohesive structures, as in  \cite{Schreiber2012},  \cite{Schreiber2013a}, \cite{Shulman2017}).   

What significance does the synthetic geometry program offer towards formalizing an approach to fundamental physics based on rewriting systems?  The limiting Wolfram model rulial multiverse has the formal structure of an  $\left(\infty, 1\right)$-topos, constructed from rewriting systems. The objects of this topos are rulial multiway systems as $\infty$-groupoids, which in turn, are formal homotopy spaces, upon which constructions relevant to physics may be  realized. Various Wolfram model combinatorial structures  such as (ordinary) multiway systems, branchial graphs, causal graphs, etc.,  are obtained as sections or restrictions of the limiting rulial multiway system, and can therefore be endowed with associated spatial structures. In this way,  pregeometric constructions in the Wolfram model can be internalized within a suitable $\infty$-topos, within which, spaces and specific constructions relevant to physics can then be realized  (such an approach is also being  pursued in  \cite{Schreiber2012},  \cite{Schreiber2013a},  \cite{schreiber2016higher} with the  aim of  formalizing quantum field theories on these toposes).  

Importantly, our work here advances the formal connection between abstract rewriting systems and higher homotopy spaces. This framework suggests a formal basis for models seeking to construct geometric structures starting from purely pregeometric notions.  
This can be potentially useful for formally justifying choices of underlying spacetime discretization adopted by contemporary models of quantum gravity  \cite{isham1995structural}, \cite{gibbs1995small}, \cite{isham2000some}, \cite{dowker2006causal}, \cite{smolin2008three}, \cite{smolin2018we}.  More specifically, and to different degrees, most models of quantum gravity presume   that a well-defined continuum limit exists, where one might retrieve smooth classical spaces \cite{rideout2009emergence}, \cite{rovelli2008loop}, \cite{loll2019quantum}. However, that may not always be as straightforward.   
 The insight that our formalization brings to this picture is that once a given discretization scheme can be expressed type-theoretically or as a multiway rewriting system, then, the homotopy structures admissible within those discretization schemes can potentially help determine spatial properties emerging from those models.

\subsection{Categorical Quantum Mechanics, Higher Homotopies and Topological Field Theories}

In previous works \cite{Gorard2020c}, \cite{Gorard2021a}  we have demonstrated how the categorical approach to quantum mechanics or CQM \cite{abramsky2009categorical}, \cite{coecke2018picturing}, and in particular, ZX calculus, can be formally expressed within the framework of our multiway rewriting systems, where the associated Wolfram model rules are precisely the equational rules of ZX calculus. The states of this multiway system consist of ZX diagrams, and edges correspond to diagram transformations. In fact, this ZX multiway system is formally an embedding space of ZX diagrams with Wolfram model evolution corresponding to double-pushout (DPO) rewriting formulated using selective adhesive categories.  Moreover, the multiway system, thus obtained, is a monoidal category with the tensor product of ZX diagrams being compatible with the tensor product of multiway states. Also,  the rulial  multiway system defined by applying all possible rules of  ZX calculus to a given ZX diagram forms a subcategory of the category of directed cospans of selective adhesive rules.

How does this picture extend in the light of higher homotopies? A potential extension here is that this leads to Higher Categorical Quantum Mechanics (HCQM), which would also serve as the starting points for Topological Quantum Field Theories (TQFTs)  \cite{baez1995higher},  \cite{lurie2009classification}.  The crucial idea here would be to formalize homotopies of the rulial multiway system as cobordisms within a suitable higher category. More of this will be reported in upcoming work.

\subsection{The Role of Homotopies as Multiway Completions for Observer-Enacted Measurement}

The Wolfram model takes the perspective that an observer has to be a part of the underlying multiway system (possibly as a subgraph spread across branches)  \cite{Wolfram2020}.  In this view, measurement is consequently  the process of the observer conflating  parallel threads of multiway history with a single evolution leading to the ``illusion" of a unique sequential thread of time. In algebraic logic and proof theory,  the process of reaching a definitive state can be operationalized via completions or mergers of branches in a proof tree \cite{eklof1971model}, \cite{ghilardi2013sheaves}. 

How exactly might completions be operationalized within Wolfram model multiway systems?  

In order to address this, let us first state the following two definitions:

\begin{definition}
A  `\textit{Completion}' in the Wolfram model is defined as an additional rule or collection of rules introduced into a multiway system that brings it closer to causal invariance (some multiway systems can be made causal invariant by adding only a finite number of completions). 
\end{definition}

In theoretical computer science, completions such as the Knuth-Bendix completion \cite{knuth1970simple}, \cite{huet1981complete}, are commonly used in automated theorem-proving algorithms, as a means of forcing confluence within equational rewriting systems. 
 
\begin{definition}
In the Wolfram model an `\textit{Observer}' is defined as any ordered sequence of non-intersecting level surfaces of a universal time function, defined over a directed acyclic graph. 
\end{definition} 

In the case of a causal network, this corresponds to part or whole of a foliation of spacetime (and therefore to an   observer embedded in a particular reference frame). In the case of a multiway evolution graph, this corresponds to part or whole of a foliation of branchial time.

In  \cite{Gorard2020a},  it was hypothesized that completion procedures such as Knuth-Bendix completions could be used to model measurement (potentially quantum measurement), allowing one to ``collapse'' superpositions of states in branchial space.  While there may be several ways of implementing the Knuth-Bendix algorithm on multiway systems,  it turns out that the explicit framework of higher homotopies in multiway systems, introduced in this work, gives a   formal  and computational procedure for implementing completions on multiway systems. In a sense, when we introduce higher homotopies between paths in a  rewriting system, we are effectively introducing a new notion of equivalence between terms in the ambient type. This is yet another notion of equivalence from the notion of state equivalence introduced by the computation rule for multiway system. These equivalences define homotopy classes of multiway states. Hence, the process of measurement (by an observer) is one which extends the type constructor for the rewriting system by introducing a kind ``extended'' computation that specifies equivalences inducing higher homotopies between paths in the multiway system. What all of this seems to suggest is that any thorough description of the (quantum)  measurement problem in models of physics based on rewriting systems would, at the very least, require inclusion of higher homotopies and a higher categorical description.

\subsection{A Homotopical Interpretation of Graph and Hypergraph Limits for Discrete Spacetime Models}

In this work we have argued that pregeometric homotopy types serve as the starting point for synthetic definitions of topology and geometry. Such a framework suggests a way for constructing geometric structures relevant to discrete models of spacetime; potentially, for those related to certain theories of quantum gravity. In particular, this provides formal criteria for justifying when different choices of underlying spacetime discretization are admissible. For instance, most contemporary approaches to quantum gravity propose various models of discrete spacetime; either by starting with a given  classical manifold, and discretizing it piecewise, or starting with a simplicial complex corresponding to a known topology. The former can be achieved either by sprinkling causal sets \cite{dowker2006causal} or discretizing an  $n$-dimensional Lorentzian manifold using causal $n$-simplices \cite{loll2019quantum}; whereas, the latter is a way to  obtain spin networks that are subsequently evolved to give spin foams \cite{rovelli2008loop}. The limiting graphs or hypergraphs of these constructions are then obtained as prescribed scaling limits which recover classical spacetime manifolds. 

There are, however, at least two issues with taking such scaling limits, that need to be addressed. The first one is of a technical nature. Namely, what if there are  local obstructions to scaling resulting in local "hairs" or other irregularities? The obvious answer would be to introduce appropriate regularization schemes. But such schemes will depend on a case-by-case basis. Hence, how should one characterize generic criteria to ensure that such procedures yield regular geometries?  The second issue is conceptual and in some sense even more pressing. If classical spaces are to truly emerge from  fundamental discrete building blocks, then shouldn't one do away with any ambient geometry upon which sprinkling or piecewise discretization is performed, or even have the luxury of using a pre-defined simplicial structure with given topology? In that case, the existence of a well-defined continuum limit is not always guaranteed. In other words, if one strictly commits to pregeometric building blocks, how should one obtain continuous spaces?   
   
This is where the homotopical rewriting systems formalized in this work become useful. So long as the process for composing discrete entities in any given spacetime model admits higher homotopical mappings as part of its construction towards an $\infty$-groupoid, Grothendieck's homotopy hypothesis posits that the resulting limiting structure is a topological space.   
Furthermore, with additional cohesivity conditions (referring to locality and connectedness) one gets smooth spaces as is done in homotopical constructions of synthetic geometry. Indeed, the synthetic geometry perspective turns the relation between topology and groupoids around, placing precedence on the latter. Sections of the limiting object one thus constructs can then be interpreted as $n$-fold groupoids, corresponding to $n$-homotopy cells. Contractions and projections of these sections yield various models of topological or geometric spaces.
In other words, graph and hypergraph limits often assumed to exist in discrete spacetime models are well-defined only when  relevant homotopical structures are admissible within those models. Therefore, constructing spaces comes down to a systematic algorithm for gluing together higher homotopies, for which, the homotopical rewriting systems explored in this work offer such a construction. The constructive algorithm employed here makes use of $n$-fold categories, using which, scaling relations on associated graphs or hypergraphs can be interpretable as scaling of hypercubes. It is important to note though, that behind any admissible scaling is its associated homotopical structure. 

\begin{figure}[ht]
\centering
\includegraphics[width=0.4\textwidth]{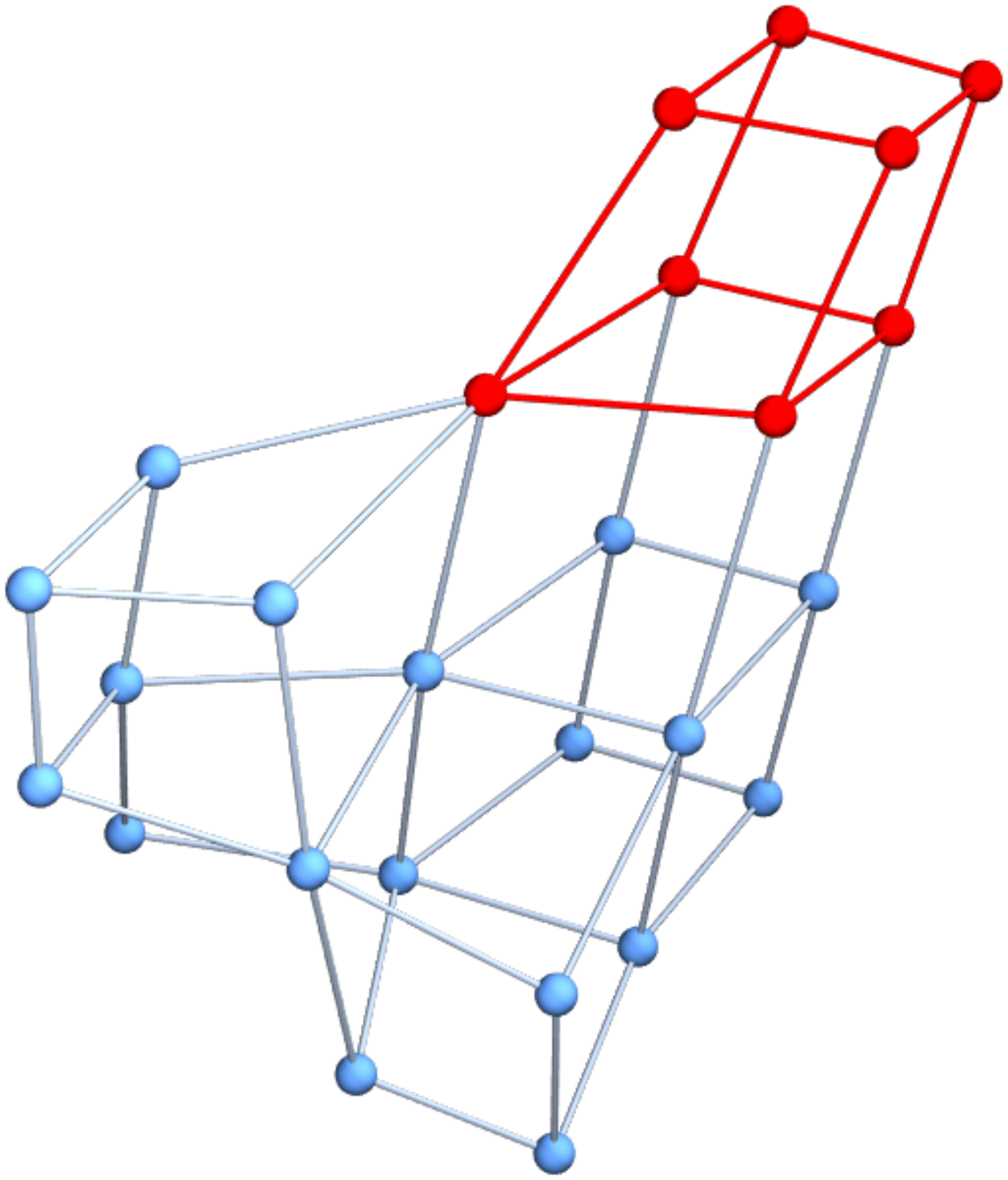}\hspace{0.1\textwidth}
\includegraphics[width=0.4\textwidth]{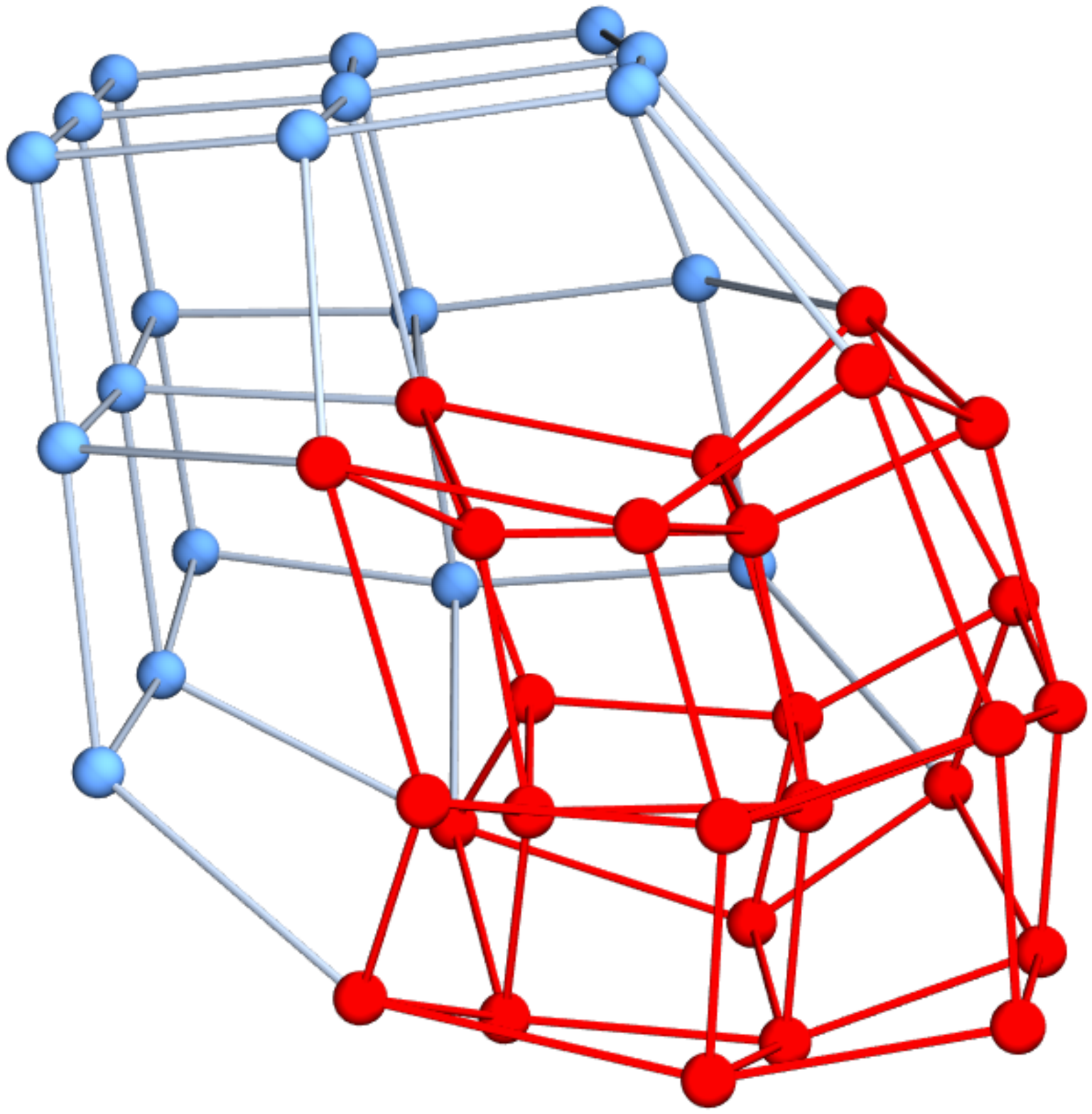}
\caption{The figure on the left side shows a rewriting system starting from an initial cubical graph that is updated via the rewriting rule that glues new cubical graphs on any of its faces. The right side shows how an initial cubical graph can be iteratively subdivided into further cubes. The red highlighting in both figures indicates  instantaneous application of their respective rewriting rules. }
\label{fig:app1}
\end{figure}

As schematic examples, let us consider the two graph rewritings shown in Figure \ref{fig:app1}. The rewriting rule for the figure on the left side involves adding a cube (with matching faces) at each update; whereas, the rule for the figure on the right involves decomposing every cube into eight "smaller" cubes. If viewed solely in terms of graph limits, the structure on the left will appear to "grow" irregularly in three dimensions and in the asymptotic limit one may be able to perform a course-graining procedure where local patches of the limiting graph can be treated as sampled approximations of three dimensional continuous maps. On the other hand, viewed in terms of a homotopical construction, the specific rewriting rule implemented here corresponds to performing compositions of 3-morphisms of a $3$-fold category. So long as higher homotopies are admissible (using additional rewriting rules) the limiting homotopical structure yields an $\infty$-groupoid, which is then amenable to synthetic constructions of topology and geometry. And only then would a section of that total space yield the three dimensional object that the graph limit procedure approximates to. Likewise, from a graph limit perspective, the construction on the right side of  Figure \ref{fig:app1} limits to something that "fills" the initial cube and approximates to a three dimensional solid cube. However, as a homotopical construction what allows for this structure are 3-morphisms, which with additional homotopies and cohesivity conditions enable continuous spaces. 

The insight one gains from this is that once a given discretization scheme can be expressed type-theoretically (as multiway systems), the homotopy structures admissible within that discretization scheme help determine spatial properties emerging from those models. The specific spaces that result still depend on particular rewriting rules used. However, the flexibility to model spaces using such systems makes them useful pregeometric tools that may potentially complement conventional spacetime constructions as those used in causal sets, causal dynamical triangulation or spin foam models. This will be elaborated further in future work.

\section{Conclusions }

The Wolfram model originates from the idea that the basic building blocks of the universe are fundamentally discrete and are governed by computational rules. Consequently, structures relevant to physics are proposed as combinatoric constructions of non-deterministic rewriting systems. A central theme of this paper was to provide a formal foundation for the Wolfram model using homotopy types. In particular, we showed an explicit representation of multiway rewriting systems using type constructors as well as an algorithm for constructing higher homotopies on these systems. Following this, all Wolfram model constructions can be represented using homotopy types.  In doing so, we have argued how spatial structures and geometry can be formally obtained from pregeometric type-theoretic constructions. This new formulation of the Wolfram model  thus provides a new computational framework that operationalizes  Wheeler's original intuition of physics from pregeometry. 
 
More specifically, in this work we have explicitly demonstrated how higher categorical constructs arise naturally in Wolfram model non-deterministic rewriting systems, the so-called  `multiway systems'.  We have shown how higher homotopies induced on multiway systems via specific rewriting rules correspond to morphisms of an $n$-fold category, and have  demonstrated an explicit computational algorithm for constructing homotopies on multiway systems. Additionally,  we have also provided the formal correspondence of these systems to $n$-fold categories and $n$-fold groupoids. Interestingly, the $n \to \infty$ limit of the rulial multiway system yields an $\infty$-groupoid, with the latter being relevant from the point of view of Grothendieck's homotopy hypothesis. Furthermore, we have shown how this construction extends to the classifying space of rulial multiway systems, which forms a multiverse of multiway systems and carries the formal structure of an ${\left(\infty, 1\right)}$-topos.  This correspondence to spaces and higher structures potentially offers  a new way to understand how the kinds of spatial structures relevant to the foundations of physics may emerge from abstract combinatorial systems.  Indeed, the pregeometric constructions elucidated here can be internalized within a suitable $\infty$-topos, within which, spaces and constructions relevant to physics can potentially be realized.  A related program seeking to formalize quantum field theories from cohesive $\infty$-toposes also borrows heavily from synthetic geometry  \cite{Schreiber2012}, \cite{Schreiber2013a}, \cite{schreiber2016higher}.  

In fact, this approach to the foundations of physics from synthetic geometry is   very much related to (and has benefited from) recent developments at the foundations of mathematics, where  results from homotopy type theory and models of  infinity-categories have been extrapolated to formalize, among other things,    definitions of higher geometric structures (in particular, cohesive topological and geometric spaces)  \cite{Schreiber2013a},  \cite{Shulman2017},      \cite{riehl20152},  \cite{Riehl2017a},  \cite{riehl2018elements}.  These constructions of spatial structures from homotopy types are precisely what the Wolfram model (as a type theory) operationalizes. And, analogous to universes of mathematics being  formalized in homotopy type theories \cite{Program2013}, \cite{Shulman2017},   the Wolfram model can be thought of as a `constructivist' approach to the foundations of physics. 

An extremely important consequence of constructing spaces and geometry synthetically is that it removes any need to make ad hoc assumptions about  geometric attributes to local constructs of a model (such as assuming a background geometric space or discretization defined using pre-assigned geometric data). Instead, geometry (in the form of cohesive structures) is inherited functorially  by  higher structures, and as a consequence, is naturally induced upon local structures by taking sections or projections of the total space. In fact, the multiway rewriting systems, considered in this work,  are abstract higher categorical objects constructed as a homotopy type. Being purely syntactic constructions, they do not hinge upon any a priori notion of space or geometry. In other words, these are  strictly pregeometric structures. It is only in well-defined limits that these pregeometric structures lead to   $\infty$-categories, where topology and geometry can be synthetically realized. Subsequently, all Wolfram model constructions obtained as sections (branchial spaces) or restrictions (causal graphs) of these limiting structures are endowed with geometry induced via the global geometry of the rulial multiway system. 

A key issue we have addressed in this work is to formally relate abstract rewriting systems to homotopy spaces. The former are purely discrete syntactic structures, whereas the latter serve as the starting point for synthetic definitions of topology and geometry. Importantly, such a framework suggests a way for constructing geometric structures relevant to physics, starting from purely pregeometric models. This can be particularly useful for formally justifying different choices of underlying spacetime discretization adopted by various  approaches to modeling  quantum gravity.  More specifically, most contemporary models of discrete spacetime assume that some well-defined continuum limit exists, where one can retrieve smooth classical spaces. However, this may not always be true. The insight that our formalization brings to this picture is that once a given discretization scheme can be expressed type-theoretically (as multiway systems), the homotopy structures admissible within those discretization schemes help determine spatial properties emerging from those models.  
 
It is interesting to note, that just as the framework of monoidal categories has led to a formal (and useful) reformulation of certain structures of quantum mechanics (in what is now called Categorical Quantum Mechanics or CQM); likewise, the framework of homotopy type theory and higher categorical structures seem the natural extension to formalizing the foundations of physics. In fact, the building blocks of CQM are combinatorial constructions, the so-called string diagrams, which represent a quantum process algebra and are formalized within dagger symmetric monoidal categories. Similarly, the building blocks of the Wolfram model are multiway rewriting systems, which represent formal models of computation or proof systems and which, within the framework of homotopy type theory, syntactically define spatial structures upon which the  foundations of physics can be built. Indeed, it has already been shown in earlier work that the Wolfram model multiway system of ZX diagrams  is a formal embedding space of ZX updating processes of CQM, with the multiway rewriting rules precisely corresponding to the equational rules of ZX calculus \cite{Gorard2020c}.  Hence, in the light of what we have demonstrated in this paper, a natural extension would be to formalize higher categorical and in particular extended topological quantum field theories  \cite{lurie2009classification}  using  rulial multiway system as cobordisms within a suitable $\infty$-category. More of this will be reported in upcoming work. 

Finally, we reckon that the framework of higher category-theoretic  combinatorial constructions developed in this work, may also be relevant to several other research programs investigating questions concerning the foundations of physics.  These include quantum mechanics on toposes, higher pre-quantum geometry, higher gauge field theories,  Lie algebroids, and branes with  higher spin excitations associated to cohomological hierarchies of the  Whitehead tower;  among others  \cite{isham2000some},  \cite{Schreiber2012},  \cite{Schreiber2013a}, \cite{schreiber2016higher}.  
In earlier work, connections to categorical quantum mechanics and causal set theory have at best only utilized ordinary category-theoretic structures of the Wolfram model \cite{Gorard2020c}, \cite{Gorard2020b}. Now,  this  new extended   formulation of the Wolfram model  based on higher categories offers new formal and computational tools to investigate higher symmetries and extended structures in topological field theories starting from homotopical multiway systems, as the new building blocks.  In this sense, the Wolfram model seeks to serve as a unifying `meta-framework' for contextualizing disparate ideas at the foundations of physics.  



%

\section*{Acknowledgments}

The authors would like to thank Stephen Wolfram for his encouragement, and for the many conversations that helped shape this work. Thanks also to Amar Hadzihasanovic, Nils A Baas, Nicolas Behr and Hatem Elshatlawy for useful feedback.

\bibliographystyle{eptcs}
\bibliography{hottrefs.bib}

\end{document}